\documentclass[mathpazo]{cicp}
\usepackage{cite}
\usepackage{bm}
\usepackage{subfigure}
\usepackage{xcolor}
\usepackage[ruled]{algorithm2e}
\usepackage{tabu}
\usepackage{multirow}

\begin{document}
\title{PFNN-2: A Domain Decomposed Penalty-Free Neural Network Method for Solving Partial Differential Equations}

\author[Sheng H L et.~al.]{Hailong Sheng\affil{1,3},
        and Chao Yang\affil{2}\comma\corrauth}
\address{\affilnum{1}\ Institute of Software,
          Chinese Academy of Sciences,
          Beijing 100190, P.R. China. \\
         \affilnum{2}\ School of Mathematical Sciences,
          Peking University,
          Beijing 100871, P.R. China. \\
         \affilnum{3}\ University of Chinese Academy of Sciences,
          Beijing 100190, P.R. China.}
\emails{{\tt hailong2019@iscas.ac.cn} (H.~Sheng),
        {\tt chao\_yang@pku.edu.cn} (C.~Yang)}

\begin{abstract}
A new penalty-free neural network method, PFNN-2, is presented for solving partial differential equations, which is a subsequent improvement of our previously proposed PFNN method \cite{sheng2021pfnn}.
PFNN-2 inherits all advantages of PFNN in handling the smoothness constraints and essential boundary conditions of self-adjoint problems with complex geometries, and extends the application to a broader range of non-self-adjoint time-dependent differential equations.
In addition, PFNN-2 introduces an overlapping domain decomposition strategy to substantially improve the training efficiency without sacrificing accuracy.
Experiments results on a series of partial differential equations are reported, which demonstrate that PFNN-2 can outperform state-of-the-art neural network methods in various aspects such as numerical accuracy, convergence speed, and parallel scalability.

\end{abstract}


\maketitle

\section{Introduction}
\label{sec1}
In recent years, neural network methods are becoming an attractive alternative for solving partial differential equations (PDEs) arising from applications such as fluid dynamics \cite{raissi2020hidden, jin2021nsfnets, hennigh2021nvidia}, quantum mechanics
\cite{han2018solving, beck2019machine, pfau2020ab}, molecular dynamics \cite{razakh2021pnd}, material sciences \cite{shukla2020physics} and geophysics \cite{li2020coupled, zhu2021general}.
In contrast to most traditional numerical approaches, methods based on neural networks are naturally meshfree and intrinsically nonlinear therefore can be applied without going through the cumbersome step of mesh generation and could be more potentially applicable to complicated nonlinear problems.
These advantages have enabled neural network methods to draw increasingly more attention
with both early studies using shallow neural networks \cite{lee1990neural, meade1994numerical, van1995neural, lagaris1998artificial, lagaris2000neural, mcfall2009artificial, mcfall2013automated} and recent works with the advent of deep learning technology \cite{raissi2019physics, karniadakis2021physics, weinan2018deep, han2018solving, long2019pde, li2020fourier, chen2020meta, lu2021learning}.

Despite of the tremendous efforts to improve the performance of neural network methods for solving PDEs, there are still issues that require further study. The first issue is related to the accuracy of neural network methods. It is still not fully understood how well a neural network can approximate the solution of a PDE in either theory or practice. Recent investigations were carried out to reveal some preliminary approximation properties of neural networks for simplified problems \cite{he2019relu, weinan2020machine, opschoor2020deep, kutyniok2021theoretical}, or to improve the accuracy of neural network methods from various aspects such as by introducing weak form loss functions to relax the smoothness constraints \cite{weinan2018deep, ming2021deep, khodayi2020varnet, kharazmi2019variational, kharazmi2020hp, zang2020weak}, by modifying the solution structure to automatically satisfy the initial-boundary conditions \cite{lagaris1998artificial, mcfall2009artificial, mcfall2013automated, liu2019solving, lyu2020enforcing}, by revising the network structure to enhance its representative capability \cite{wang2020understanding, jagtap2020adaptive, gao2021phygeonet, gao2021super}, and by introducing advanced sampling strategies to reduce the statistical error \cite{nabian2021efficient, wight2020solving}.
These improvements so far are usually effective in respective situations, but are not general enough to adapt with different types of PDEs defined on complex geometries.

Another difficulty for solving PDEs with neural network methods is related to the efficiency. It is widely known that training a neural network is usually much more costly than solving a linear system resulted from a traditional discretization of the PDE \cite{strikwerda2004finite, bathe2006finite, eymard2000finite}.
To improve the efficiency, it is natural and has been extensively considered \cite{hennigh2021nvidia, mcclenny2021tensordiffeq} to utilize distributed training techniques \cite{ben2019demystifying, farkas2020parallel} based on data or model parallelism, by which the training task is split into a number of small sub-tasks according to the partitioned datasets or model parameters so that multiple processors can be exploited.
Although this distributed training is general and successful in handling many machine learning tasks \cite{chang2018distributed, chen2019distributed, zhang2019distributed, wu2020collaborate}, it is not the most effective choice for solving PDEs because no specific knowledge of the original problem is adopted throughout the training process.

Inspired by the idea of classical domain decomposition methods \cite{smith1996pe}, it was proposed to introduce domain decomposition strategies into neural network based PDE solvers \cite{li2020deep, mercier2021coarse, li2019d3m, jagtap2020conservative, jagtap2020extended, shukla2021parallel, hu2021extended, jagtap2022physics, de2022error}
by dividing the learning task into training a series of sub-networks related to solutions on sub-domains.
This tends to be more natural than the plain distributed training approaches because by introducing domain decomposition the training of each sub-network only requires a small part of dataset related to the corresponding sub-domain, therefore can significantly decrease the computational cost.
However, these methods may still suffer from issues related to low convergence speed and poor parallel scalability due in large part to the straightforward treatment of artificial sub-domain boundaries.
In this paper, we present PFNN-2, a new penalty-free neural network method that is a subsequent improvement of our previously proposed PFNN method \cite{sheng2021pfnn}. Inheriting all advantages of PFNN in handling the smoothness constraints and essential boundary conditions of self-adjoint PDEs, PFNN-2 extends the application to a broader range of non-self-adjoint time-dependent problems and introduces a domain decomposition strategy to improve the training efficiency. In particular, PFNN-2
employs a Galerkin variational principle to transform the non-self-adjoint time-dependent equation into a weak form,
and adopts compactly supported test functions that are convenient to build on complex geometries. 
To further reduce the difficulty of the training process and improve the accuracy, PFNN-2 separates the penalty terms related to the initial-boundary constraints from the loss function, with a dedicated neural network learning the initial and essential boundary conditions, and the true solution on the rest of the domain approximated by another network that is totally disentangled from the former one.
On top of the above, PFNN-2 introduces an overlapping domain decomposition strategy to speedup multi-processor training without sacrificing accuracy, 
thanks to the penalty-free treatment of initial-boundary conditions of the decomposed sub-problems and the weak form loss function with reduced smoothness requirement.
We carry out a series of numerical experiments to demonstrate that PFNN-2 is advantageous over existing state-of-the-art algorithms in terms of both numerical accuracy and computational efficiency.

The remainder of the paper is organized as follows. In Section 2 we present the basic algorithms employed in PFNN-2 to handle the smoothness constraints and initial-boundary conditions. Following that the domain decomposition strategy adopted by PFNN-2 for parallel computing is introduced in Section 3. We then discuss and compare PFNN-2 with several related state-of-the-art neural network methods in Section 4. Numerical experiment results on a series of linear and nonlinear, time-independent and time-dependent, initial-boundary value problems on simple and complex geometries are reported in Section 5. The paper is concluded in Section 6.

\section{The PFNN-2 method}
\label{sec2}
For the ease of discussion, we first consider the case of a single sub-domain $\Omega\subset \mathbb{R}^{d_{\bm{x}}}$ ($d_{\bm{x}}\in\mathbb{N}_{+}$)
on which we try to solve the following initial-boundary value problem:
\begin{equation}
  \label{ibvp}
  \left \{
    \begin{array}{r l l}
      \dfrac{\partial u}{\partial t}
      - \nabla \cdot
        \left( \mathcal{A} \nabla u -
               \mathcal{B} \right)
      + \mathcal{C} &= 0,
      & \mbox{in}\ \Omega\times(0,T], \\[2.5mm]

      u &= r_D,\ \
      & \mbox{on}\ \Gamma_D\times(0,T], \\[2mm]

      \left( \mathcal{A} \nabla u \right)
      \cdot\bm{n} &= r_N,
      & \mbox{on}\ \Gamma_N\times(0,T], \\[2mm]

      u &= r_0,
      & \mbox{on}\ \overline{\Omega}\times\{0\},
    \end{array}
  \right.
\end{equation}
where
$\mathcal{A}$, $\mathcal{B}$ and $\mathcal{C}$ are dependent on the solution $u$ as well as the spatial and temporal coordinates $\bm{x}$ and $t$, and
$\bm{n}$ is the outward unit normal.
Both Dirichlet and Neumann boundary conditions are considered, with the corresponding boundaries
$\Gamma_D\cup\Gamma_N = \partial\Omega$ and
$\Gamma_D\cap\Gamma_N = \varnothing$.

To reduce the smoothness requirement, we introduce the following weak form
\begin{equation}
  \label{vp}
  \displaystyle \int_0^T \int_{\Omega}
  \left(
    \mathcal{A} \nabla u \cdot \nabla v +
    \Big( \frac{\partial u}{\partial t} +
          \nabla\cdot\mathcal{B} +
          \mathcal{C} \Big) v
  \right) d\bm{x}dt
  - \displaystyle \int_0^T \int_{\Gamma_N}
  r_N v d\bm{x}dt = 0,
  \quad \forall v\in\mathcal{V}
\end{equation}
with the solution $u$ belonging to a trial (hypothesis) space $\mathcal{H}$ and satisfying initial-boundary conditions
\begin{equation*}
  u = r_D,\ \mbox{on}\
  \Gamma_D\times(0,T]
  \quad\mbox{and}\quad
  u = r_0,\ \mbox{on}\
  \overline{\Omega}\times\{0\},
\end{equation*}
where $\mathcal{V}$ is a certain test space, in which any function $v$ satisfies $v=0$ on $\Gamma_D\times(0,T]$.

In PFNN-2, the approximate solution of problem (\ref{vp}) is found within the hypothesis space $\mathcal{H}$ that is formed by neural networks.
Unlike most of the existing methods \cite{raissi2019physics, karniadakis2021physics, weinan2018deep, ming2021deep, khodayi2020varnet, kharazmi2019variational, zang2020weak} that merely apply a single network to construct $\mathcal{H}$, PFNN-2 adopts two neural networks, with one network $g_{\bm{\theta}_1}$ learning the solution with the initial condition and essential boundary condition more accurately and swiftly, and another network $f_{\bm{\theta}_2}$ approximating the solution on the rest part of the domain, where $\bm{\theta}_1$, $\bm{\theta}_2$ donate the sets of weights and biases forming the two networks, respectively.
To eliminate the influence of $f_{\bm{\theta}_2}$ on 
the initial condition and essential boundary condition,
a length factor function is introduced to impose appropriate restriction, which satisfies
\begin{equation}
  \label{func_l_value}
  \left \{
    \begin{array}{l l}
      \ell=0, & \mbox{on}\
      (\overline{\Omega}\times\{0\}) \cup
      (\Gamma_D\times(0,T]), \\[2mm]
      \ell>0, & \mbox{otherwise}.
    \end{array}
  \right.
\end{equation}
With the help of the length factor function, any function in the hypothesis space $\mathcal{H}$ is constructed as follows:
\begin{equation}
  \label{aprx_solu}
  w_{\bm{\theta}}:= g_{\bm{\theta}_1} + \ell f_{\bm{\theta}_2},
\end{equation}
where $\bm{\theta}= \{\bm{\theta}_1,\bm{\theta}_2\}$.

The length factor function $\ell$ is established by utilizing spline functions, which are more flexible to adapt to various complicated geometries than the analytical functions proposed in previous works \cite{lagaris1998artificial, liu2019solving, lyu2020enforcing}.
To be specific, we divide the essential boundary $\Gamma_D$ into several segments $\{\gamma_j\}_{j=1}^{n_D}$, $n_D\in\mathbb{N}_{+}$.
For each segment $\gamma_j$, select another segment of the boundary as its companion $\gamma_{j_o}$ (note that $\gamma_{j_o}$ should not be adjacent to $\gamma_j$). Then, for each $\gamma_j$, a spline function $l_j$ is built, which satisfies
\begin{equation}
  \label{l_basis}
  \left \{
    \begin{array}{l l}
      l_j(\bm{x}) = 0, & \bm{x}\in\gamma_{j},   \\ [2mm]
      l_j(\bm{x}) = 1, & \bm{x}\in\gamma_{j_o}, \\ [2mm]
      0<l_j(\bm{x})<1, & \mbox{otherwise}.
    \end{array}
  \right.
\end{equation}
The specific details to construct $l_j$ can be found in our previous work of PFNN \cite{sheng2021pfnn}. Utilizing all the $l_j$ functions, the length factor function $\ell$, which is dependent on both space and time, is established by:
\begin{equation}
  \label{func_l}
  \ell(\bm{x},t):=
  \frac{\widetilde{\ell}(\bm{x},t)}
  {\max\limits_{(\widehat{\bm{x}},\widehat{t})
                \in \Omega\times[0,T]}
                \widetilde{\ell}(\widehat{\bm{x}},\widehat{t})},
  \quad\mbox{where}\quad
  \widetilde{\ell}(\bm{x},t) =
  \frac{t}{T}
  \prod_{j=1}^{n_D} 1-(1-l_j(\bm{x}))^{\mu},
\end{equation}
where $\mu>0$ is adopted for adjusting the shape of $\ell$. In order to avoid that the value of $\ell$ drops dramatically with the growth of $n_D$, we suggest that $\mu=n_D$.

\begin{figure}[!htb]
  \centering
  \subfigure[boundary]
  {\includegraphics[height=0.15\textheight]
  {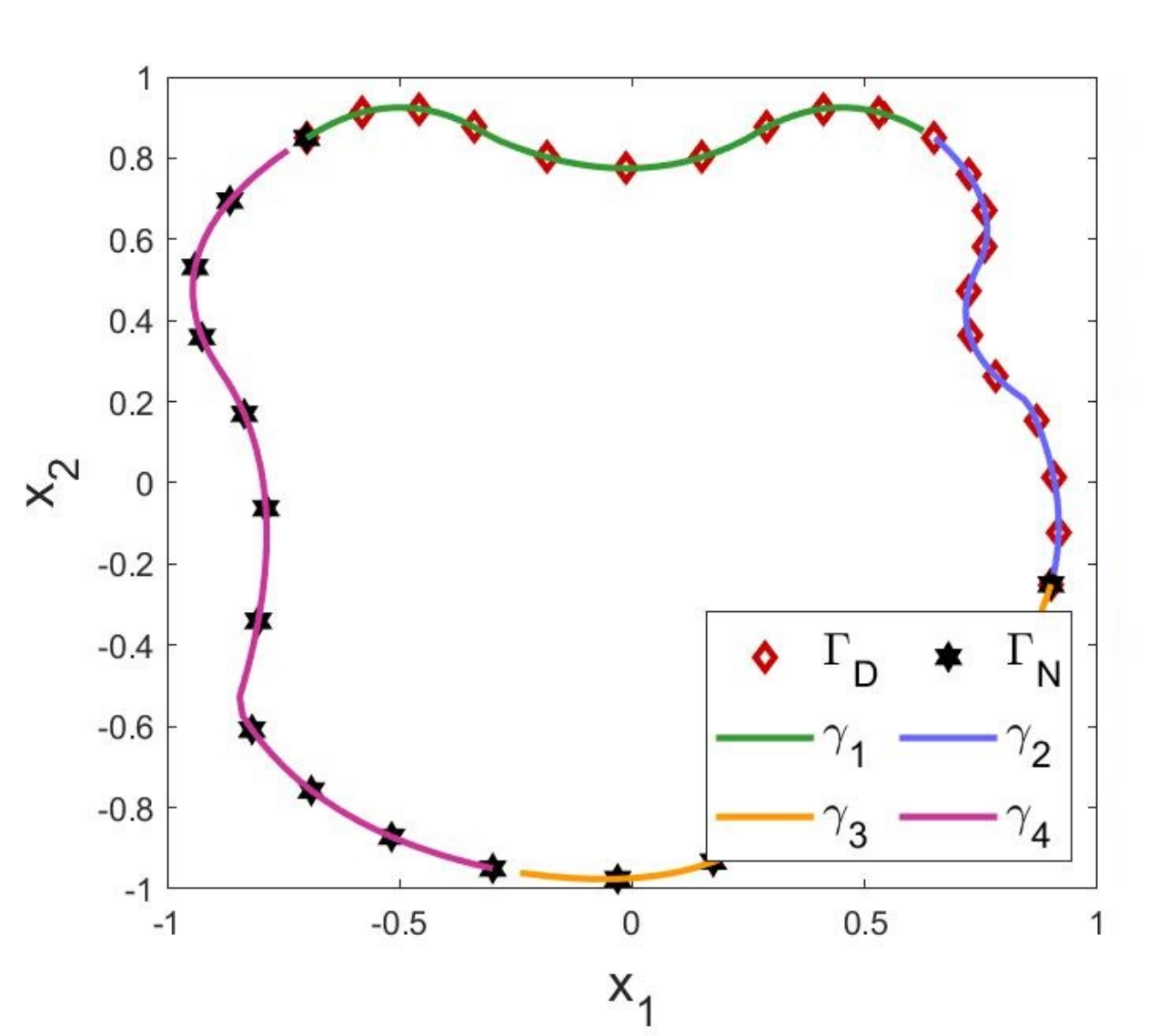}}
  \subfigure[$l_1$]
  {\includegraphics[height=0.15\textheight]
  {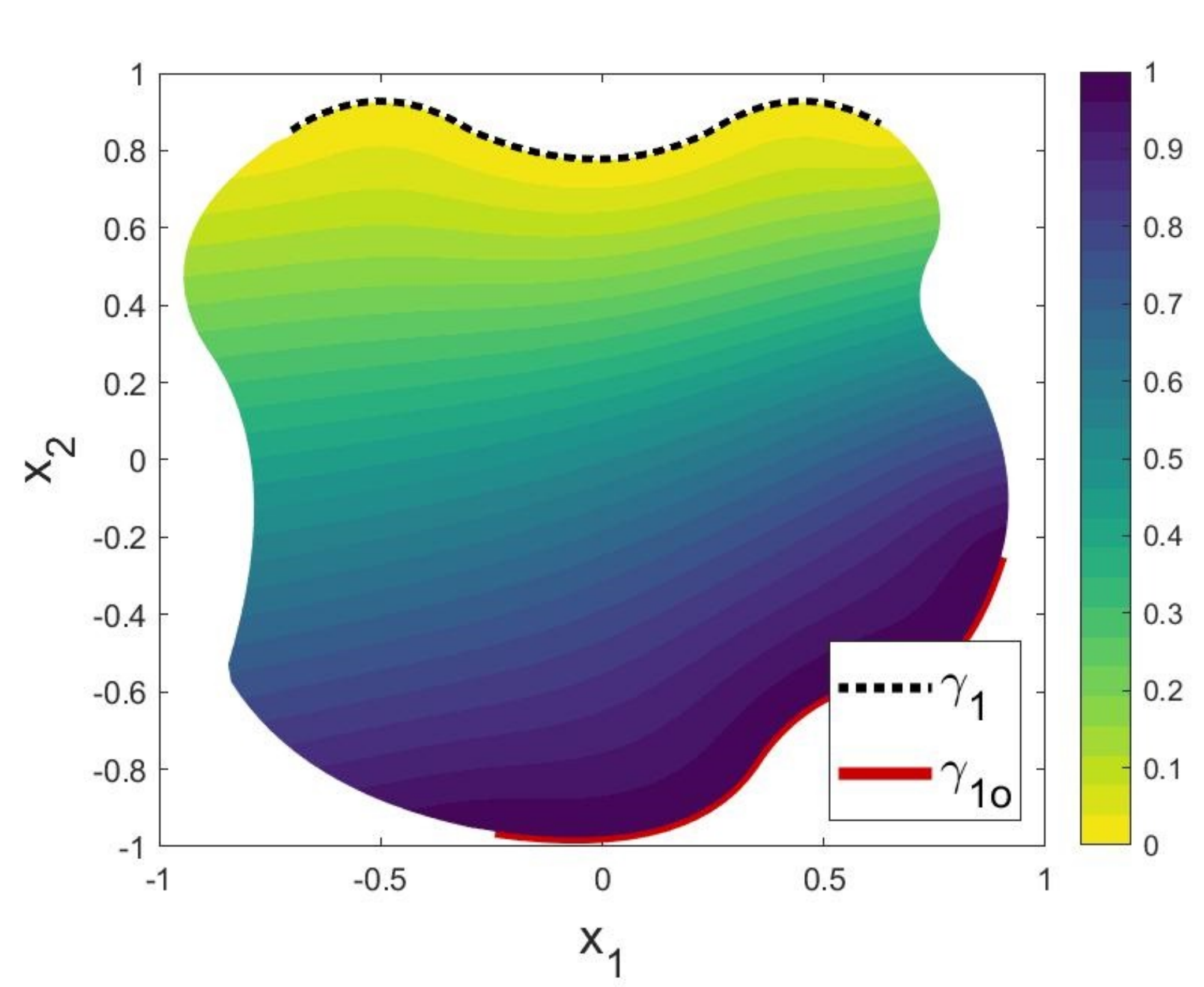}}
  \subfigure[$l_2$]
  {\includegraphics[height=0.15\textheight]
  {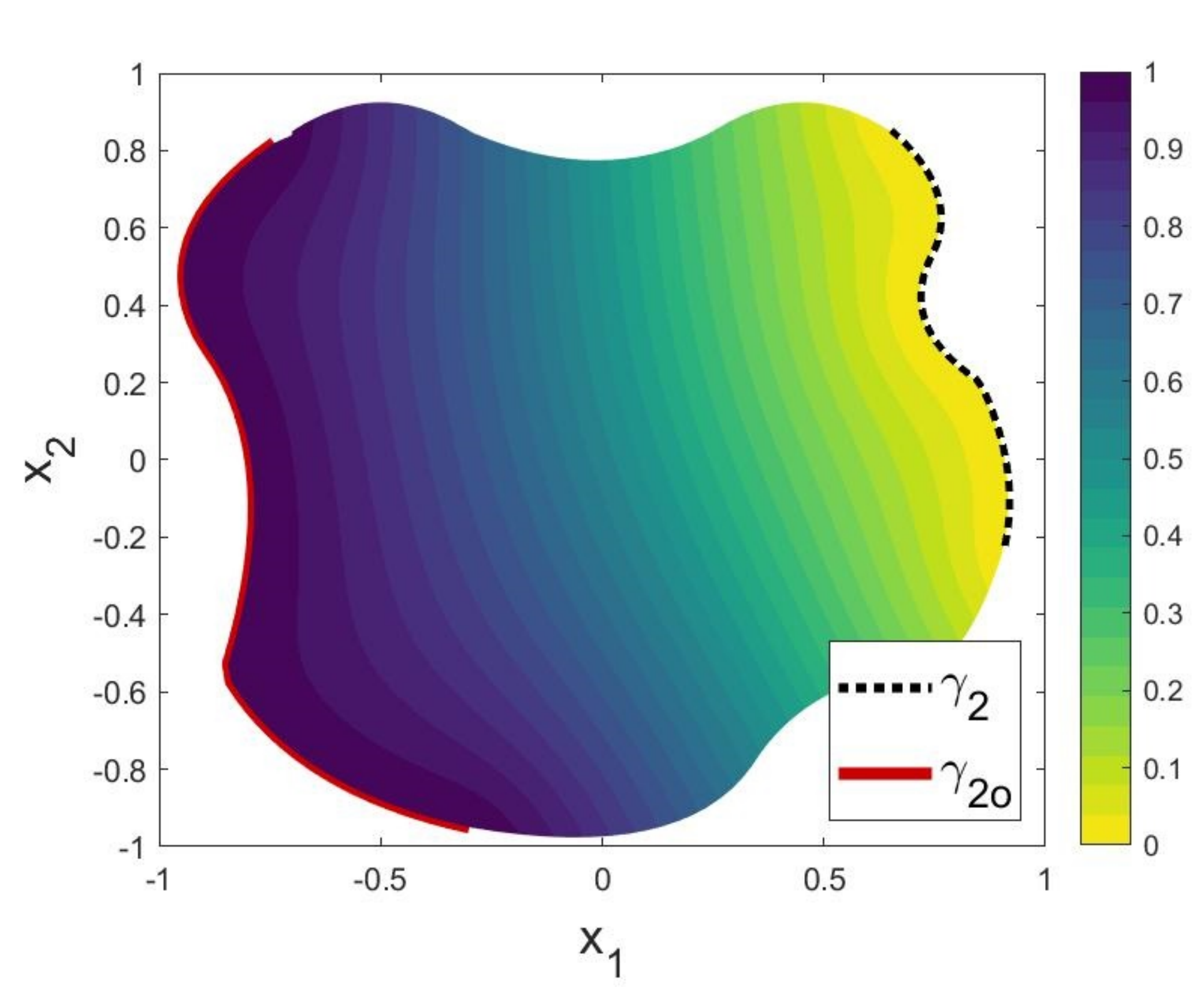}}
  \subfigure[$\ell$($\mu=2$)]
  {\includegraphics[height=0.15\textheight]
  {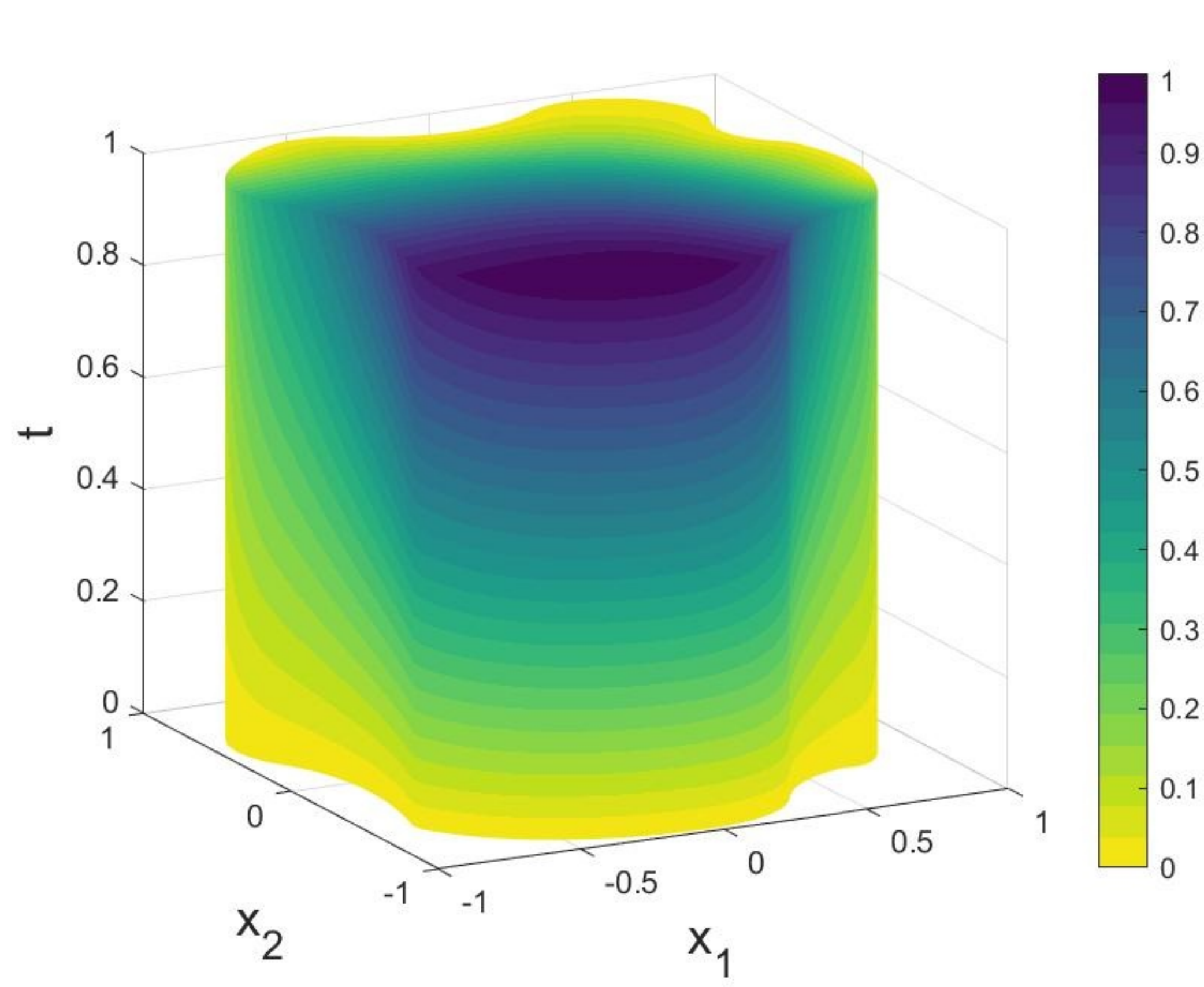}}
  \caption{An demonstration on the construction of the length factor function.}
  \label{f_length_factor}
\end{figure}

Figure \ref{f_length_factor} is an example to demonstrate how to construct the function $\ell$, in which the essential boundary $\Gamma_D$ is divided into two parts, $\gamma_1$ and $\gamma_2$. Correspondingly, another two segments, $\gamma_3$ and $\gamma_4$, are taken from the Neumann boundary $\Gamma_N$ as the companions $\gamma_{1_o}$ and $\gamma_{2_o}$, respectively. Then two spline functions $l_1$ and $l_2$ satisfying (\ref{l_basis}) can be built according to the method in PFNN \cite{sheng2021pfnn} and the length factor function is established according to (\ref{func_l}).



To improve the flexibility to complex geometries and reduce the cost to evaluation residual (\ref{vp}),
rather than the commonly utilized orthogonal polynomials \cite{kharazmi2019variational, kharazmi2020hp}, a set of compactly supported test functions
\begin{equation}
  \label{test_func}
  \small
  v_s(\bm{x},t; \widehat{\bm{x}}_{s},\widehat{t}_{s},h_s)
  := \max\left\{ 1 - \frac{|t - \widehat{t}_{s}|}
                 {h_s}, 0 \right\}
     \prod_{j=1}^{d_{\bm{x}}}
     \max\left\{ 1 - \frac{|x_j - \widehat{x}_{s,j}|}
                 {h_s}, 0 \right\},
  \ s=1,2,\cdots,n_v,\ n_v\in\mathbb{N}_{+}
\end{equation}
are employed to span the test space $\mathcal{V}$,
where $(\widehat{\bm{x}}_{s},\widehat{t}_{s})$ and $h_s$ are the center and radius of the support
$G_s = \cap_{j=1}^{d_{\bm{x}}}
 \{(\bm{x},t)| |x_j-\widehat{x}_{s,j}| < h_s \}
 \cap \{(\bm{x},t)| |t-\widehat{t}_s| < h_s \}$
of $v_s$, respectively, and
$x_j$ and $\widehat{x}_{s,j}$ are the $j$th element of $\bm{x}$ and $\widehat{\bm{x}}_s$, respectively.
The center $(\widehat{\bm{x}}_{s},\widehat{t}_{s})$ is sampled on $(\Omega\cup\Gamma_N)\times[0,T]$
while the radius $h_s$ is small to ensure that the residual can be just evaluated on a tiny region instead of the whole domain. To achieve this, the radius $h_s$ is set to:
\begin{equation}
  \label{test_func_radius}
  h_s:=
  \min\{\delta(\widehat{\bm{x}}_{s}, \Gamma_D),
        \widehat{t}_s, T-\widehat{t}_s, \overline{h} \},
\end{equation}
where $\delta(\widehat{\bm{x}}_{s}, \Gamma_D)$ is the minimum distance from $\widehat{\bm{x}}_s$ to $\Gamma_D$ so that $v_s=0$ on $\Gamma_D$, and $\overline{h}$ is the upper bound of $h_s$, which is defined as
\begin{equation*}
  \overline{h}
  = \left(\frac{|\Omega|T}{n_v}
    \right)^{\frac{1}{d_{\bm{x}}+1}}.
\end{equation*}
Such configuration keeps $\sum_{s=1}^{n_v} |G_s| \approx 2^{d_{\bm{x}}+1}|\Omega|T$ and guarantees that the supports of all the test functions could cover most part of the computing domain.

Since the two neural networks $g_{\bm{\theta}_1}$ and $f_{\bm{\theta}_2}$ are decoupled by the length factor function $\ell$, they can be trained separately.
First, $g_{\bm{\theta}_1}$ is trained to learn the initial and essential boundary condition via minimizing the loss function
\begin{equation}
  \label{loss_func_g}
  \Psi[g_{\bm{\theta}_1}]
  := \sum\limits_{(\bm{x},t)\in
     S(\Gamma_D\times(0,T])}
     \big( g_{\bm{\theta}_1}(\bm{x},t) - r_D(\bm{x},t) \big)^2
   + \sum\limits_{(\bm{x},t)\in
     S(\overline{\Omega}_i\times\{0\})}
     \big( g_{\bm{\theta}_1}(\bm{x},t) - r_0(\bm{x},t) \big)^2,
\end{equation}
where $S(\Box)$ represents a set of points sampled on $\Box$.
Then, let $g_{\bm{\theta}_1^{\ast}}$ be the obtained network after minimizing (\ref{loss_func_g}) and $\widetilde{w}_{\bm{\theta}_2} =  g_{\bm{\theta}_1^{\ast}} + \ell f_{\bm{\theta}_2}$.
An approximate solution $u^{\ast}$ of problem (\ref{vp}) can be obtained by minimizing $L$ as follows:
\begin{equation}
  \label{loss_func_f}
  u^{\ast}
  = \arg\min\limits_{\widetilde{w}_{\bm{\theta}_2}}
    L[\widetilde{w}_{\bm{\theta}_2}]
  = \arg\min\limits_{\widetilde{w}_{\bm{\theta}_2}}
    \sum_{s=1}^{n_v} R^2[\widetilde{w}_{\bm{\theta}_2}; v_s],
\end{equation}
where
\begin{equation}
  \small
  \label{residual}
  \begin{array}{r l}
    R[\widetilde{w}_{\bm{\theta}_2}; v_s]
    & := \displaystyle\frac{|\Omega| T}
       {\#S(\Omega\times (0,T])}
       \sum\limits_{(\bm{x},t)\in S(\Omega\times [0,T])}
       \Big[\mathcal{A} \nabla \widetilde{w}_{\bm{\theta}_2}
            \cdot \nabla v_s +
            \Big(\dfrac{\partial \widetilde{w}_{\bm{\theta}_2}}{\partial t} +
            \nabla\cdot\mathcal{B} +
            \mathcal{C}
            \Big) v_s
       \Big] \Big(\bm{x},t\Big) \\
     &\ - \displaystyle \frac{|\Gamma_N| T}
       {\#S(\Gamma_N\times(0,T])}
       \sum\limits_{(\bm{x},t)\in S(\Gamma_N\times[0,T])}
       \big[r_N v_s\big] \big(\bm{x},t\big)
  \end{array}
\end{equation}
is the estimate of the residual of the governing equation (\ref{vp}), $\#S(\Box)$ donates the size of $S(\Box)$.
Eq. (\ref{loss_func_f}) is taken as the loss function for training $f_{\bm{\theta}_2}$.
It can be seen that neither of the loss functions (\ref{loss_func_g}) and (\ref{loss_func_f}) involves any penalty term, therefore is easier to minimize than those with penalty terms frequently found in various neural network methods for solving partial differential equations. Besides, the training of each network only needs data samples on part of the computing domain, which has less computational cost than training a single network with data samples from the whole domain.

\section{Domain decomposition for PFNN-2}
\begin{figure}[!htb]
  \centering
  \includegraphics[width=0.9\textwidth]
  {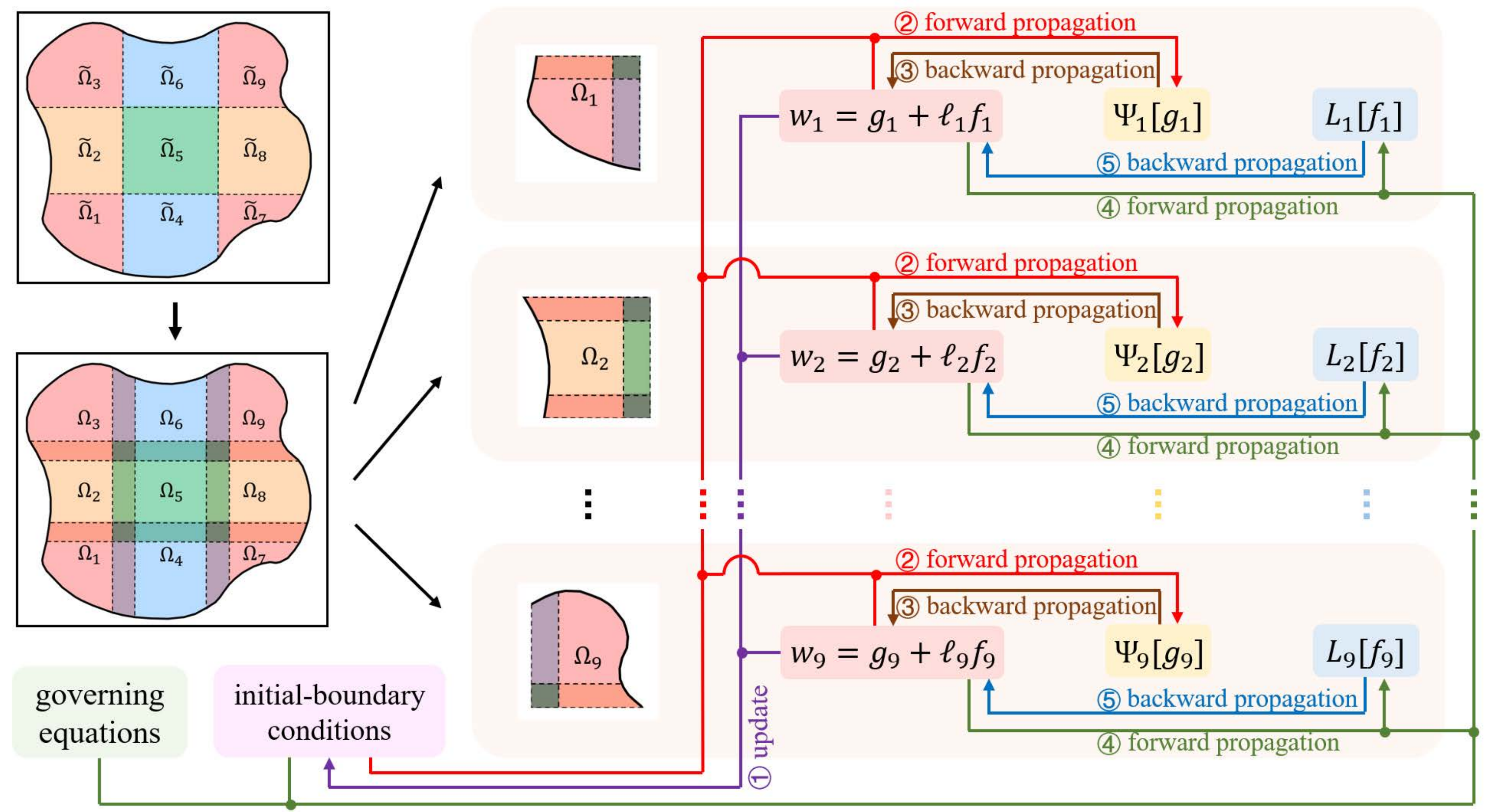}
  \caption{The domain partitioning, network structures and training process of PFNN-2.
  The computational domain is decomposed into several overlapping sub-domains and
  each sub-domain comes with a combination of two sub-networks for approximating the relative local solution.
  Each iteration of the training process is comprised of steps \textcircled{1}-\textcircled{5}.}
  \label{f_domain_decomposition}
\end{figure}

To enable parallel computing capability in PFNN-2, we first partition the computing domain into a group of non-overlapping sub-domains $\{\widetilde{\Omega}_i\}_{i=1}^m$, and then extend each sub-domain to obtain a series of overlapping sub-domains $\{\Omega_i\}_{i=1}^m$,
where $m\in\mathbb{N}_{+}$ is the number of sub-domains.
Figure \ref{f_domain_decomposition} shows an example of the overlapping domain decomposition partition for the case of
$m=9$.
Starting with a suitable initial guess $u^0$,
problem (\ref{ibvp}) can be solved by constructing a solution sequence $u^k$ that converges to the true solution, where $k=1,2,\cdots$ represents the count of iteration. And the approximate solution for each iteration is formed by solving a series of sub-domain problems concurrently on parallel computers.
In particular, at the $k$-th iteration the sub-domain problem defined on $\Omega_i$ is
\begin{equation}
  \label{ibvp_sub}
  \begin{array}{l}
    \left \{
    \begin{array}{r l l}
      \dfrac{\partial u_i^k}{\partial t}
      - \nabla \cdot
        \left( \mathcal{A} \nabla u_i^k -
               \mathcal{B} \right)
      + \mathcal{C}
      & = 0,
      & \mbox{in}\ \Omega_i\times(0,T], \\[2mm]

      u_i^k
      & = u^{k-1},
      & \mbox{on}\ \Gamma_{i}\times(0,T], \\[2mm]

      u_i^k
      & = r_D,
      & \mbox{on}\ \Gamma_{i,D}\times(0,T], \\[2mm]

      \left( \mathcal{A} \nabla u_i^k \right)
      \cdot\bm{n}
      & = r_N,
      & \mbox{on}\ \Gamma_{i,N}\times(0,T], \\[2mm]

      u_i^k
      & = r_0,
      & \mbox{on}\ \overline{\Omega}_i\times\{0\}
    \end{array}
    \right. 
  \end{array}
\end{equation}
where $\Gamma_{i,D}= \partial\Omega_i \cap \Gamma_D$,
$\Gamma_{i,N}= \partial\Omega_i \cap \Gamma_N$, and
$\Gamma_{i}= \partial\Omega_i \backslash (\Gamma_D \cup \Gamma_N)$.
After solving the sub-domain problems for iteration $k$, the approximate solution defined on the whole domain is composed as $u^{k}\big|_{\widetilde{\Omega}_i} = u_i^{k}$.




To construct the solution sequence, instead of solving the sub-domain problem (\ref{ibvp_sub}) directly, PFNN-2 finds an approximate solution to the corresponding sub-domain problem in the weak form
\begin{equation}
  \label{vp_sub}
  \displaystyle \int_0^T \int_{\Omega_i}
  \left(
    \mathcal{A} \nabla u_i^k \cdot \nabla v +
    \Big( \frac{\partial u_i^k}{\partial t} +
          \nabla\cdot\mathcal{B} +
          \mathcal{C} \Big) v
  \right) d\bm{x}dt
  - \displaystyle \int_0^T \int_{\Gamma_{i,N}}
  r_N v d\bm{x}dt = 0,
  \quad \forall v\in\mathcal{V}_i,
\end{equation}
where $u_i^k$ belongs to a hypothesis space $\mathcal{H}_i$ and satisfies the initial-boundary conditions
\begin{equation*}
  \label{ib_con_sub}
  \begin{array}{l}
    \left \{
    \begin{array}{r l l}

      u_i^k
      & = u^{k-1},
      & \mbox{on}\ \Gamma_{i}\times(0,T], \\[2mm]

      u_i^k
      & = r_D,
      & \mbox{on}\ \Gamma_{i,D}\times(0,T], \\[2mm]

      u_i^k
      & = r_0,
      & \mbox{on}\ \overline{\Omega}_i\times\{0\},
    \end{array}
    \right. 
  \end{array}
\end{equation*}
and $\mathcal{V}_i$ is a certain test space.

Since the weak form (\ref{vp_sub}) is similar to (\ref{vp}), we can utilize the method presented in the previous section to establish the hypothesis space $\mathcal{H}_i$ and test space $\mathcal{V}_i$, 
for all $i=1,2,\cdots,m$.
Specifically on sub-domain $\Omega_i$, $\mathcal{V}_i$ is spanned by a set of compactly supported functions $\{v_{i,s}\}_{s=1}^{n_v^i}$ ($n_v^i\in\mathbb{N}_{+}$),
and any function $w_i\in\mathcal{H}_i$ is constructed according to
\begin{equation}
  \label{aprx_sub}
  w_i = g_i + \ell_i f_i,
\end{equation}
where $g_i$ and $f_i$ are two sub-networks that
play roles similar to $g_{\bm{\theta}_1}$ and $f_{\bm{\theta}_2}$,
and $\ell_i$ is a length factor function that decouples the two sub-networks.
Note that here and hereafter, we drop the subscripts representing the learnable parameters for brevity.

Starting from a proper initial state $w_i^0 = g_i^0 + \ell_i f_i^0$, each local approximate solution $w_i$ is updated by a dedicated processor to produce a sequence $w_i^k$ converging to the corresponding local true solution, where $g_i^0$ and $f_i^0$ are the initial states of sub-networks $g_i$ and $f_i$, $k=1,2,\cdots$ is the count of iteration.
In the $k$-th iteration, $w_i^k$ is acquired by solving the corresponding sub-problem (\ref{vp_sub}).
And following that, all the $w_i^k$ are merged together to form the approximate solution $w^k$ on the whole domain according to $w^k\big|_{\widetilde{\Omega}_i} = w_i^k$.
To illustrate the solving process accurately and intuitively, Algorithm 1 and Figure \ref{f_domain_decomposition} are provided, which show that each iteration is comprised of two stages:
\begin{itemize}
  \item {\bf Online stage (line 4-13 of Algorithm 1):} As demonstrated in step \textcircled{1} of Figure \ref{f_domain_decomposition}, each processor conducts communication to the immediate neighbors to acquire the current interior boundary condition,
  where $w_i$ satisfies $w_i^k = w^{k-1} = w_j^{k-1}$, on $(\Gamma_i\cap \Omega_j) \times (0,T]$ $(j\neq i)$.

  \item {\bf Offline stage (line 15-23 of Algorithm 1):} First, as shown in steps \textcircled{2}-\textcircled{3} of Figure \ref{f_domain_decomposition}, each sub-network $g_i$ is trained on the corresponding processor to learn the local solution by minimizing the loss function
      \begin{equation}
        \label{loss_func_sub_g}
        \small
        \begin{array}{l}
          \Psi_i^k[g_i]
          := \displaystyle\sum\limits_{(\bm{x},t)\in
             S(\Gamma_i\times(0,T])}
             \big( g_i(\bm{x},t) - w^{k-1}(\bm{x},t) \big)^2
           + \displaystyle\sum\limits_{(\bm{x},t)\in
             S(\Gamma_{i,D}\times(0,T])}
             \big( g_i(\bm{x},t) - r_D(\bm{x},t) \big)^2 \\[6mm]
           \qquad\quad
           + \displaystyle\sum\limits_{(\bm{x},t)\in
             S(\overline{\Omega}_i\times\{0\})}
             \big( g_i(\bm{x},t) - r_0(\bm{x},t) \big)^2.
        \end{array}
      \end{equation}
      Let $g_i^k$ be the acquired sub-network after minimizing (\ref{loss_func_sub_g}) and $\widetilde{w}_i^k=g_i^k + \ell_i f_i$.
      Then, as shown in steps \textcircled{4}-\textcircled{5} of Figure \ref{f_domain_decomposition}, each sub-network $f_i$ is trained on the corresponding processor to approximate the local solution on the remainder of the sub-domain via minimizing
      \begin{equation}
        \label{loss_func_sub_f}
        \small
        \begin{array}{l}
          L_i^k[f_i]:= \displaystyle\sum_{s=1}^{n_v^i}
          \Bigg[
            \displaystyle\frac{|\Omega_i| T}
            {\#S(\Omega_i\times (0,T])}
            \sum\limits_{(\bm{x},t)\in S(\Omega_i\times [0,T])}
            \Big[\mathcal{A} \nabla \widetilde{w}_i^k
                 \cdot \nabla v_{i,s} +
                 \Big(\dfrac{\partial \widetilde{w}_i^k}{\partial t} +
                 \nabla\cdot\mathcal{B} +
                 \mathcal{C}
                 \Big) v_{i,s}
            \Big] \Big(\bm{x},t\Big) \\[2mm]
            \qquad\qquad\qquad\
            - \displaystyle \frac{|\Gamma_{i,N}| T}
            {\#S(\Gamma_{i,N}\times(0,T])}
            \sum\limits_{(\bm{x},t)\in S(\Gamma_{i,N}\times[0,T])}
            \big[r_N v_{i,s}\big] \big(\bm{x},t\big)
          \Bigg]^2.
        \end{array}
      \end{equation}
      Afterwards, the obtained sub-network $f_i^k$ is utilized to form the local approximate solution $w_i^k = g_i^k + \ell_i f_i^k$.
      And finally, the approximate solution on the whole domain is composed as $w^k\big|_{\widetilde{\Omega}_i} = w_i^k$.
\end{itemize}

\begin{algorithm}[!htb]
\caption{Computing the $i$-th local approximate solution in PFNN-2. }
\LinesNumbered
\small
\KwIn{$i$, $m$, $g_i^0$ $f_i^0$, $\ell_i$, $K_o$ (number of outer iteration), $K_i$ (number of inner iteration).}
\KwOut{$w_i^{K_o}$.}
Let $w_i^0 = g_i^0 + \ell_i f_i^0$.  \\
\For{$k=1,2,\cdots,K_o$}
{
    /* Online stage: line 4 to 13 */ \\
    {
    \For{$j=1,2,\cdots,m$}
    {
        \If{$i<j$}
        {
            Send $w^{k-1}\big|_{\Gamma_j \cap \Omega_i} = w_i^{k-1}$
            to processor $j$. \\
            Receive $w^{k-1}\big|_{\Gamma_i \cap \Omega_j} = w_j^{k-1}$
            from processor $j$.
        }
        \If{$i>j$}
        {
            Receive $w^{k-1}\big|_{\Gamma_i \cap \Omega_j} = w_j^{k-1}$
            from processor $j$.\\
            Send $w^{k-1}\big|_{\Gamma_j \cap \Omega_i} = w_i^{k-1}$
            to processor $j$.
        }
    }
    }
    /* Offline stage: line 15 to 23 */ \\
    {
    \For{$t=1,2,\cdots,K_i$}
    {
        Forward propagation: evaluate the loss function $\Psi_i^k[g_i]$ according to (\ref{loss_func_sub_g}). \\
        Backward propagation: calculate the gradients of $\Psi_i^k[g_i]$ and minimize $\Psi_i^k[g_i]$ to get the updated $g_i^k$. \\
    }
    \For{$t=1,2,\cdots,K_i$}
    {
        Forward propagation: evaluate the loss function $L_i^k[f_i]$ according to (\ref{loss_func_sub_f}). \\
        Backward propagation: calculate the gradients of $L_i^k[f_i]$ and minimize $L_i^k[f_i]$ to get the updated $f_i^k$. \\
    }
    Let $w_i^k = g_i^k + \ell_i f_i^k$.
    }
}
\end{algorithm}

\section{Discussion and comparison with other related methods}
In this section, we review several state-of-the-art neural network methods for solving PDEs and compare them with PFNN-2. For simplicity we consider a model boundary value problem
\begin{equation}
  \label{bvp}
  \left \{
    \begin{array}{r l}
      -\Delta u = r_{\Omega},
      & \mbox{in}\ \Omega, \\[2mm]
      u = r_{\Gamma},
      & \mbox{on}\ \Gamma:=\partial\Omega.
    \end{array}
  \right.
\end{equation}

The most commonly applied neural network algorithms are based on the idea of least squares; examples can be found in e.g., Refs. \cite{van1995neural, lagaris1998artificial, mcfall2009artificial, raissi2019physics, raissi2020hidden, sirignano2018dgm}. In a least-squares based method, an approximate solution of (\ref{bvp}) is found by minimizing the loss function in the least-squares form as follows:
\begin{equation}
  \label{loss_least_squares}
  L[w]
  := \dfrac{1}{\#S(\Omega)}
     \sum\limits_{\bm{x}\in S(\Omega)}
     \big( -\Delta w(\bm{x}) - r_{\Omega}(\bm{x})\big)^2
   + \dfrac{\beta}{\#S(\Gamma)}
     \sum\limits_{\bm{x}\in S(\Gamma)}
     \big( w(\bm{x}) - r_{\Gamma}(\bm{x})\big)^2,
\end{equation}
where $\beta>0$ is a penalty coefficient.
Despite of being straightforward to implement, this approach may encounter  issues of low convergence speed and poor accuracy, primarily because the complicated loss function usually requires simultaneously approximating the high-order derivatives of the solution and the additional initial and boundary conditions, thus imposing extreme difficulty for optimization.
In addition, the training of the highly nonlinear neural network may also demand large computational cost and lead to low computational efficiency.

\subsection{Accuracy improvement}
Several neural network methods were proposed to enhance the accuracy by lowering the smoothness requirement so that high-order differential operators involved in the loss function are eliminated.
These methods can be roughly classified into two categories, involving weak forms based on the Ritz and Galerkin variational principles, respectively.

In approaches based on the Ritz variational principle such as the ones presented in Refs. \cite{weinan2018deep, ming2021deep, wang2020mesh, huang2021augmented},  the original problem in the strong form (\ref{bvp}) is transformed to an energy minimization problem
\begin{equation}
  \label{vp_ritz}
  u = \arg\min_{w\in\mathcal{H}}
  \displaystyle \int_{\Omega}
  \left( \dfrac{1}{2} |\nabla w|^2 -
         r_{\Omega} w \right) d\bm{x},
\end{equation}
where $\mathcal{H}$ is a certain hypothesis space built by neural networks.
In this way, the high-order derivatives are avoided in the loss function.
However, the scope of applications of these methods is usually relatively narrow, confined to self-adjoint elliptic problems.

In methods based on the Galerkin variational principle, an approximate solution is found within a hypothesis space $\mathcal{H}$ satisfying
\begin{equation}
  \label{vp_galerkin}
  \displaystyle \int_{\Omega}
  \big( \nabla u \cdot \nabla v -
        r_{\Omega} v \big) d\bm{x} = 0,
  \quad \forall v \in\mathcal{V},
\end{equation}
where $\mathcal{V}$ is a certain test function space.
Although having a wider applications than Ritz-type approaches, existing Galerkin-based neural network methods are usually in lack of flexibility and efficiency. 
For example, orthogonal polynomials can be employed as the test functions \cite{kharazmi2019variational, kharazmi2020hp}, but the methods are only applicable to simple geometries for which orthogonal polynomials are known.
It is also possible to take compactly supported piecewise linear functions to construct the test space with low cost \cite{khodayi2020varnet},
which however requires a grid to arrange the locations of the test functions and is again hard to adapt to complex geometries.
In some recent works \cite{zang2020weak, bao2020numerical}, an adversarial network was introduced to evaluate the maximum weighted loss in an adaptive manner, but the method is often hard to converge because of the difficulty in balancing the efforts devoted to training the generative network for the approximate solution and the adversarial network for the test function.

Apart from relaxing the smoothness constraints, there are another series of neural network methods that aim to improve the accuracy by removing the boundary constraints \cite{lagaris1998artificial, mcfall2009artificial, mcfall2013automated, liu2019solving, lyu2020enforcing}.
The basic idea is to establish the approximate solution in the form of
\begin{equation}
  \label{aprx_solu_ana}
  w:= \mathcal{G} + \mathcal{L} f,
\end{equation}
where $\mathcal{G}$ is manufactured to satisfy boundary conditions, $f$ is a network for learning the true solution within the domain, and $\mathcal{L}$ is introduced to eliminate the influence of the network $f$ on the boundary.
The main disadvantage of these methods is lacking a universal and flexible strategy to build the functions $\mathcal{G}$ and $\mathcal{L}$.
For handling a specific problem on a particular domain, most of the works \cite{lagaris1998artificial, liu2019solving, lyu2020enforcing} rely on a dedicated function in the analytical form, which usually cannot be applied to any other case.
Although some efforts have been made \cite{mcfall2009artificial, mcfall2013automated} by introducing spline functions so that the method can be applied to a broader range of computing domains, the flexibility is still very limited.

The proposed PFNN-2 method integrates the advantages of all the aforementioned methods. 
With the introduction of the weak form loss function, the burdensome task of approximating the high-order derivatives are avoided.
By separating the task of learning initial and essential boundary conditions from other training tasks, the constraints are dealt with in a more effective manner, and the subsequent training process is eased because the penalty terms in the corresponding loss function are eliminated.
The structures of the approximate solution comprised of neural networks and spline functions as well as the compactly supported test functions have a high flexibility, which can extend the applicability of the method to a wide range of PDEs.

\subsection{Domain decomposition}
In a classical numerical method for solving PDEs, the idea of domain decomposition can be introduced by transforming the solution of the original problem into solving a number of simpler sub-problems defined on a series of sub-domains $\{\Omega_i\}_{i=1}^m$ ($m\in\mathbb{N}_{+}$).
Starting from an appropriate initial guess $u^0$, a solution sequence $u^k$ ($k=0,1,2,\cdots$) is built to approach to the true solution of the original problem.
In particular, for problem (\ref{bvp}), at the $k$-th iteration, $u^k$ is formed by composing a set of local approximate solutions $\{u_i^k\}_{i=1}^m$, which are acquired by solving sub-domain problems
\begin{equation}
  \label{bvp_sub}
  \left \{
  \begin{array}{l l}
    \mathfrak{D}_i(u_i^k) = 0,
    & \mbox{in}\ \Omega_i, \\[2mm]
    \mathfrak{B}_i(u_i^k; u^{k-1}) = 0,
    & \mbox{on}\ \Gamma_i:=\partial\Omega_i,
  \end{array}
  \right.
  \qquad i=1,2,\cdots,m
\end{equation}
in parallel, where $\mathfrak{D}_i$ and $\mathfrak{B}_i$ are certain differential operators, whose forms depend on the specific method employed.

Previously proposed domain decomposition methods for neural network based solution of PDEs can be categorized into two major
classes: overlapping and non-overlapping methods. Figure \ref{f_domain_decomposition} has already illustrated the difference between the two. Among existing works, DeepDDM \cite{li2020deep} and D3M \cite{li2019d3m} are overlapping methods, and cPINNs \cite{jagtap2020conservative, shukla2021parallel} and XPINNs \cite{jagtap2020extended, shukla2021parallel} are  non-overlapping ones. Depending on the partition strategies adopted, as well as other characteristics, the differential operators $\mathfrak{D}_i$ and $\mathfrak{B}_i$ in these methods are in various forms, as briefly listed below.
\begin{itemize}
  \item DeepDDM:
      \begin{equation}
        \label{bvp_sub_deepddm}
        \mathfrak{D}_i(u_i^k):=
        -\Delta u_i^k - r_{\Omega}
        \quad\mbox{and}\quad
        \mathfrak{B}_i(u_i^k):=
        \left\{
          \begin{array}{l l}
            u_i^k - r_{\Gamma}, & \mbox{on}\
            \Gamma_i\cap\Gamma, \\[2mm]
            u_i^k - u^{k-1},    & \mbox{on}\
            \Gamma_i\backslash\Gamma.
          \end{array}
        \right.
      \end{equation}

  \item D3M:
      \begin{equation}
        \label{bvp_sub_d3m}
        \mathfrak{D}_i(u_i^k):=
        \left[
          \begin{array}{l}
            \bm{\tau}_i^k + \nabla u_i^k \\[2mm]
            \nabla\cdot\bm{\tau}_i^k - r_{\Omega}
          \end{array}
        \right]
        \quad\mbox{and}\quad
        \mathfrak{B}_i(u_i^k):=
        \left\{
          \begin{array}{l l}
            u_i^k - r_{\Gamma}, & \mbox{on}\
            \Gamma_i\cap\Gamma, \\[2mm]
            u_i^k - u^{k-1},  & \mbox{on}\
            \Gamma_i\backslash\Gamma.
          \end{array}
        \right.
      \end{equation}
      where $\bm{\tau}_i^k$ is an additional variable to approximate the gradient of the solution $u_i^k$.

  \item cPINNs:
      \begin{equation}
        \label{bvp_sub_cpinns}
        \begin{array}{l}
          \mathfrak{D}_i(u_i^k):=
          -\Delta u_i^k - r_{\Omega}
          \,\ \mbox{and}\,\

          \mathfrak{B}_i(u_i^k):=
          \left\{
          \begin{array}{l l}
            u_i^k - r_{\Gamma}, & \mbox{on}\
            \Gamma_i\cap\Gamma, \\[2mm]
            \left[
            \begin{array}{l}
              u_i^k - u^{k-1} \\[2mm]
              \nabla(u_i^k - u^{k-1}) \cdot\bm{n}
            \end{array}
            \right],
            & \mbox{on}\ \Gamma_i\backslash\Gamma.
          \end{array}
          \right.
        \end{array}
      \end{equation}

  \item XPINNs:
      \begin{equation}
        \label{bvp_sub_xpinns}
        \begin{array}{l}
          \mathfrak{D}_i(u_i^k):=
          -\Delta u_i^k - r_{\Omega}
          \,\ \mbox{and}\,\

          \mathfrak{B}_i(u_i^k):=
          \left\{
          \begin{array}{l l}
            u_i^k - r_{\Gamma}, & \mbox{on}\
            \Gamma_i\cap\Gamma, \\[2mm]
            \left[
            \begin{array}{l}
              u_i^k - u^{k-1} \\[2mm]
              \Delta u_i^k  - \Delta u^{k-1}
            \end{array}
            \right],
            & \mbox{on}\ \Gamma_i\backslash\Gamma.
          \end{array}
          \right.
        \end{array}
      \end{equation}
\end{itemize}

In general, the aforementioned methods solve sub-problems (\ref{bvp_sub}) by introducing $m$ sub-networks $\{w_i\}_{i=1}^m$ to approximate the local solutions on $\{\Omega_i\}_{i=1}^m$, respectively.
After initialized properly to form the initial states $\{w_i^0\}_{i=1}^m$, the $m$ sub-networks are trained in parallel to generate solution sequences $\{w_i^k\}_{i=1}^m$ converging to $m$ local true solutions, where $k=1,2,\cdots$ is the iteration count, and
$w_i^k$ is acquired by minimizing the loss function
\begin{equation}
  \label{loss_sub}
    L_i^k[w_i]
    := \dfrac{1}{\#S(\Omega_i)}
       \sum\limits_{\bm{x}\in S(\Omega_i)}
       \big\| \left[\mathfrak{D}_i(w_i)\right]
              (\bm{x}) \big\|_2^2
     + \dfrac{\beta}{\#S(\Gamma_i)}
       \sum\limits_{\bm{x}\in S(\Gamma_i)}
       \big\| \left[\mathfrak{B}_i(w_i; w^{k-1})\right]
              (\bm{x}) \big\|_2^2.
\end{equation}
Correspondingly, a solution sequence $w^k$ defined on the whole domain is formed by combining all the updated sub-networks $\{w_i^k\}_{i=1}^m$ together.


In contrast to the existing domain decomposition based neural network methods that utilize loss functions in strong form (\ref{loss_sub}), PFNN-2 employs weak form loss functions to avoid the explicit computation of high order derivatives of the approximate solution, and therefore has a higher computing efficiency. In addition to that, thanks to the penalty-free property, PFNN-2 is able to reduce the difficulty of the network training process and improve the sustained accuracy.
As demonstrated in later numerical experiments, all these advantages make PFNN-2 a more suitable candidate for the case of multiple sub-domains, especially when the number of sub-domains is large.

\section{Numerical Experiments}
In this section, we carry out a series of numerical experiments to investigate the performance of the proposed PFNN-2 method, and compare it with several existing state-of-the-art approaches.
Four test cases are utilized to examine the accuracy and efficiency of different numerical methods.
In all tests involving neural networks, the ResNet model \cite{he2016deep} with sinusoid activation function is utilized for constructing the network, in which each ResNet block is comprised of two fully connected layers and a residual connection. Note that PFNN-2 utilizes two neural networks instead of one.
To make fair comparison between various methods, we try to adjust the number of ResNet blocks and the width of hidden layers so that the total degree of freedoms (i.e., parameters including weights and biases for neural network methods, and unknowns for traditional methods) is kept to a similar level.
The L-BFGS optimizer \cite{nocedal1980updating} is applied for training the neural networks. All the experiments are conducted on a computer equipped with Intel Xeon Gold 6130 CPUs. And all calculations are done with 32-bit floating point computations. The relative $\ell_2$-error is taken to measure the accuracy of the approximate solutions. All the results obtained with neural network methods are reported based on ten times of independent experiments.


\subsection{Anisotropic convection-diffusion equation on an L-shaped domain}
The first test case is an anisotropic convection-diffusion equation
\begin{equation}
  \label{pde_acd}
  -\nabla \cdot \mathcal{A}\nabla u + \mathcal{B}\cdot\nabla u + \mathcal{C} = 0
\end{equation}
defined on an L-shaped domain $[-2,2]^2\backslash [0,2]^2$, where
\begin{equation*}
  \mathcal{A} = \left[
  \begin{array}{c c}
    -(x_1+x_2)^2 + 16 & x_1^2 - x_2^2 \\[2mm]
    -x_1^2 + x_2^2    & -(x_1-x_2)^2 + 16
  \end{array}
  \right],
  \qquad
  \mathcal{B} = \left[
  \begin{array}{c}
    (x_1 - x_1)^3 \\[2mm]
    (x_1 + x_2)^3
  \end{array}
  \right].
\end{equation*}
We set the solution to be $ u(r,\theta) = r(r-1)(r-2) \sin(2\theta/3-\pi/3)$, where $(r,\theta)$ is the polar coordinate, and let
 $\mathcal{C}$ be derived accordingly.
We apply Neumann boundary condition on $\Gamma_N = (\{-2\}\times[-2,2]) \cup ([-2,2]\times\{-2\})$ and Dirichlet boundary condition on the rest of the boundary $\Gamma_D = \partial\Omega\backslash \Gamma_N$.
The main difficulty of solving this problem is the singularity near the inner corner of the L-shaped domain.

We run the test with only one processor to examine the performance of various methods without domain decomposition.
It should be noted that not all neural network methods based on variational principles can be applied to this problem.
In particular, some methods such as Deep Ritz \cite{weinan2018deep} and Deep Nitsche \cite{ming2021deep} is only suitable for self-adjoint problems, some such as VPINNs \cite{kharazmi2019variational} and hp-VPINNs  \cite{kharazmi2020hp} require orthogonal polynomials to form the test spaces, which are hard to build on irregular domains.
Therefore, in this test, we only examine methods including the traditional linear finite element method, the least-squares neural network method, a variational neural network method called VarNet \cite{khodayi2020varnet} (which takes compactly supported piecewise linear functions as test functions) and PFNN-2.
For PFNN-2, we sample 1200 points in $\Omega$ and 80 points on $\Gamma_N$ as the centers to form the test functions according to the strategy specified in section 2, and sample 160 points on $\Gamma_D$. For the other neural network methods, the training sets with similar size are adopted.
The maximum number of iterations for training neural networks is set to 6000 epochs.

\begin{table}[!htb]
  \caption{Accuracy of various methods for solving the anisotropic convection-diffusion equation on an L-shaped domain.}
  \label{t_linear_acd_results}
  \scalebox{0.78}{
  \renewcommand{\arraystretch}{1.1}
  \begin{tabu}{c|c|c|c|c}
    \tabucline[1pt]{-}
    Method    & finite element & least-squares & VarNet & PFNN-2 \\\tabucline[1pt]{-}
    Network   & -              & 3 blocks (width:15) & 3 blocks (width:15) & 2$\times$3 blocks (width:10) \\
    \hline
    DOFs      & 1242           & 1261          & 1261   & 1182 \\
    \hline
    \multirow{3}*{\shortstack{Relative\\[1.5mm] error}}
              & 2.32e-03       & $\,\,\,\,\,\,\beta=10$ \ 1.38e-02$\pm$6.12e-03 & $\,\,\,\,\,\,\beta=10$ \ 1.56e-02$\pm$3.43e-03 & \textbf{9.18e-04$\pm$3.37e-04
} \\
              &                & $\,\,\,\beta=100$ \ 6.64e-03$\pm$3.02e-03 & $\,\,\,\beta=100$ \ 6.71e-03$\pm$1.38e-03 \\
              &                & $\,\beta=1000$ \ 1.06e-02$\pm$4.93e-03 & $\,\beta=1000$ \ 6.82e-03$\pm$1.98e-03 \\
    \tabucline[1pt]{-}
  \end{tabu}}
\end{table}

\begin{figure}[!htb]
  \centering
  \subfigure[finite element]
  {\includegraphics[width=0.24\textwidth]
  {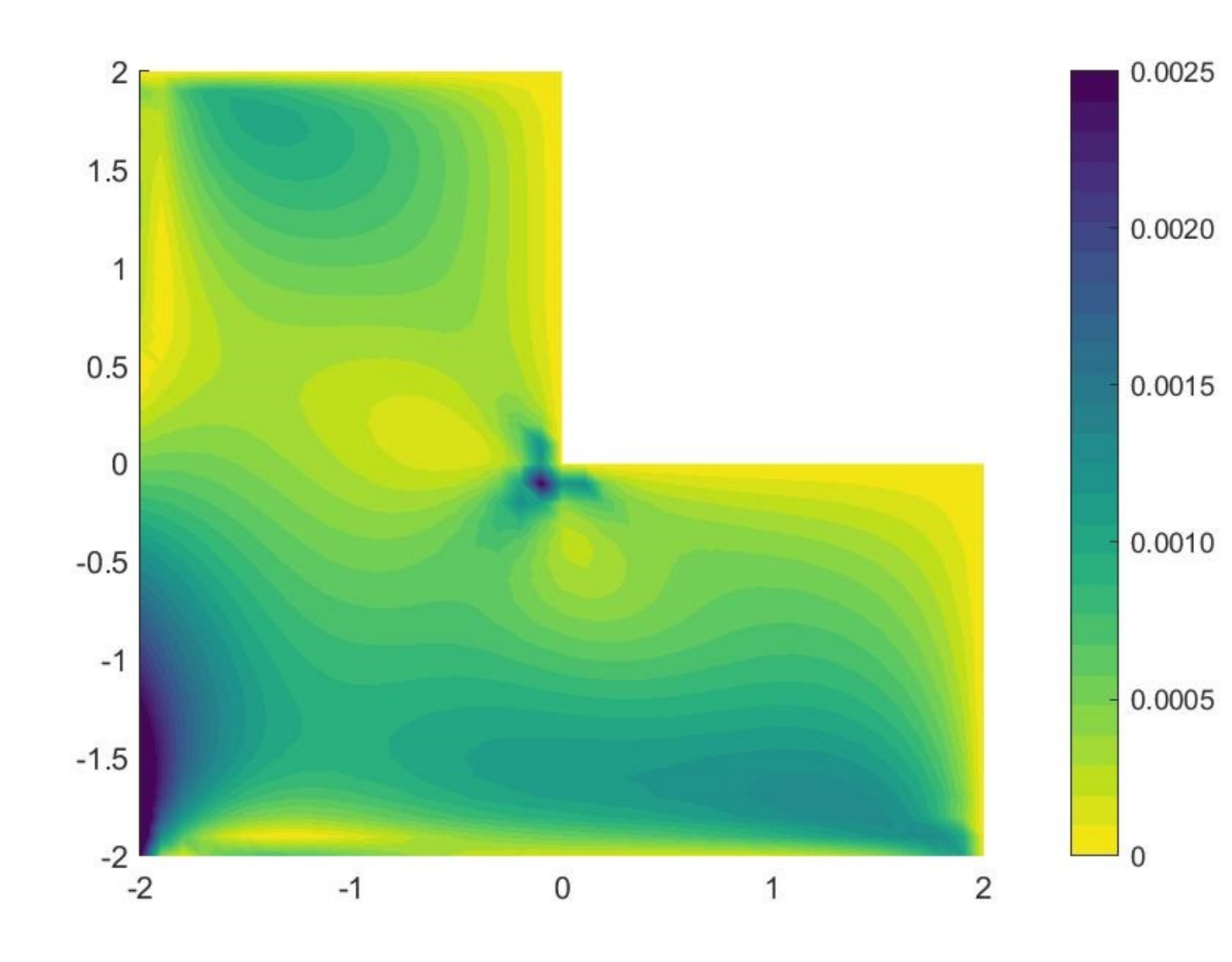}}
  \subfigure[least-squares ($\beta$=100)]
  {\includegraphics[width=0.24\textwidth]
  {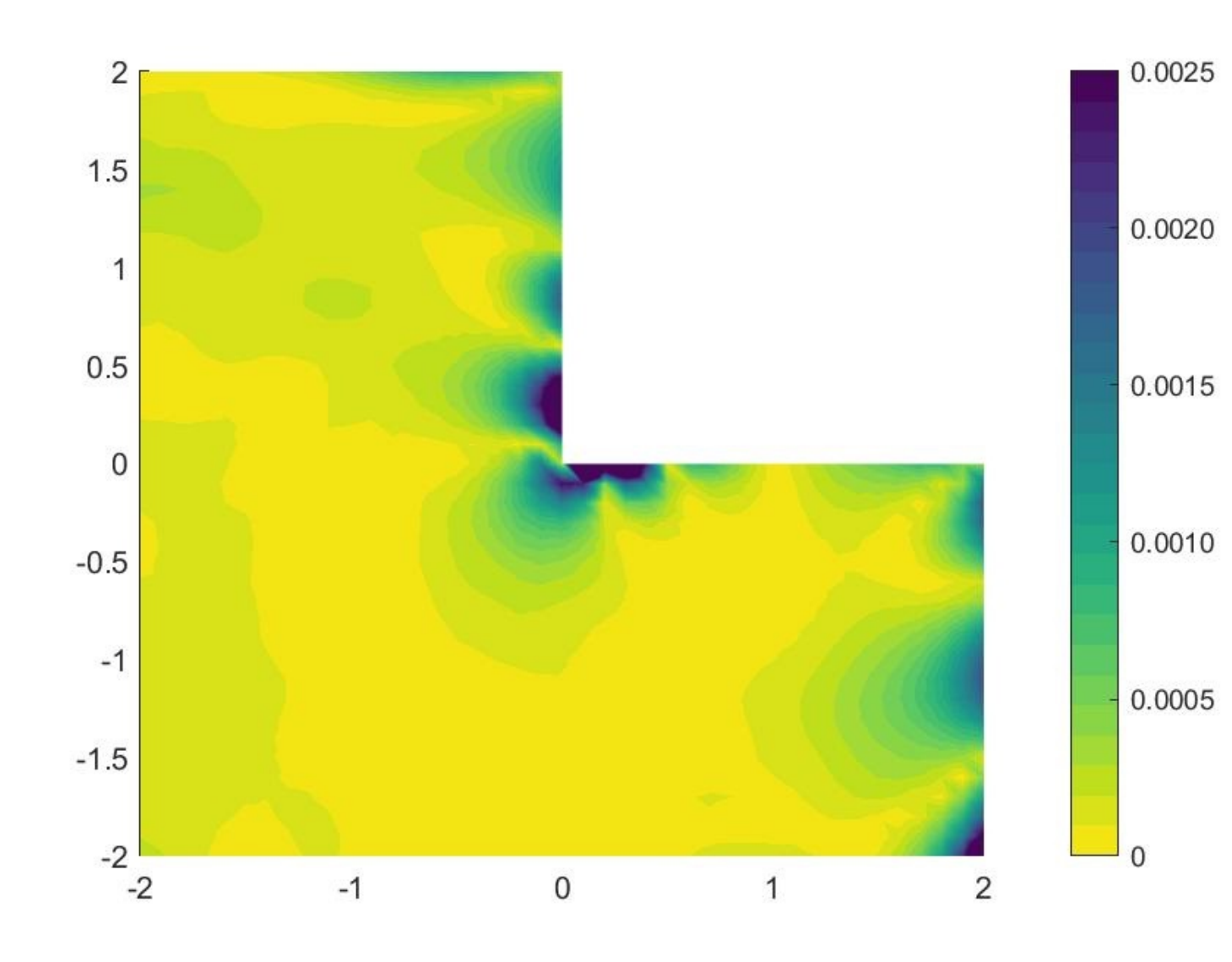}}
  \subfigure[VarNet ($\beta$=100)]
  {\includegraphics[width=0.24\textwidth]
  {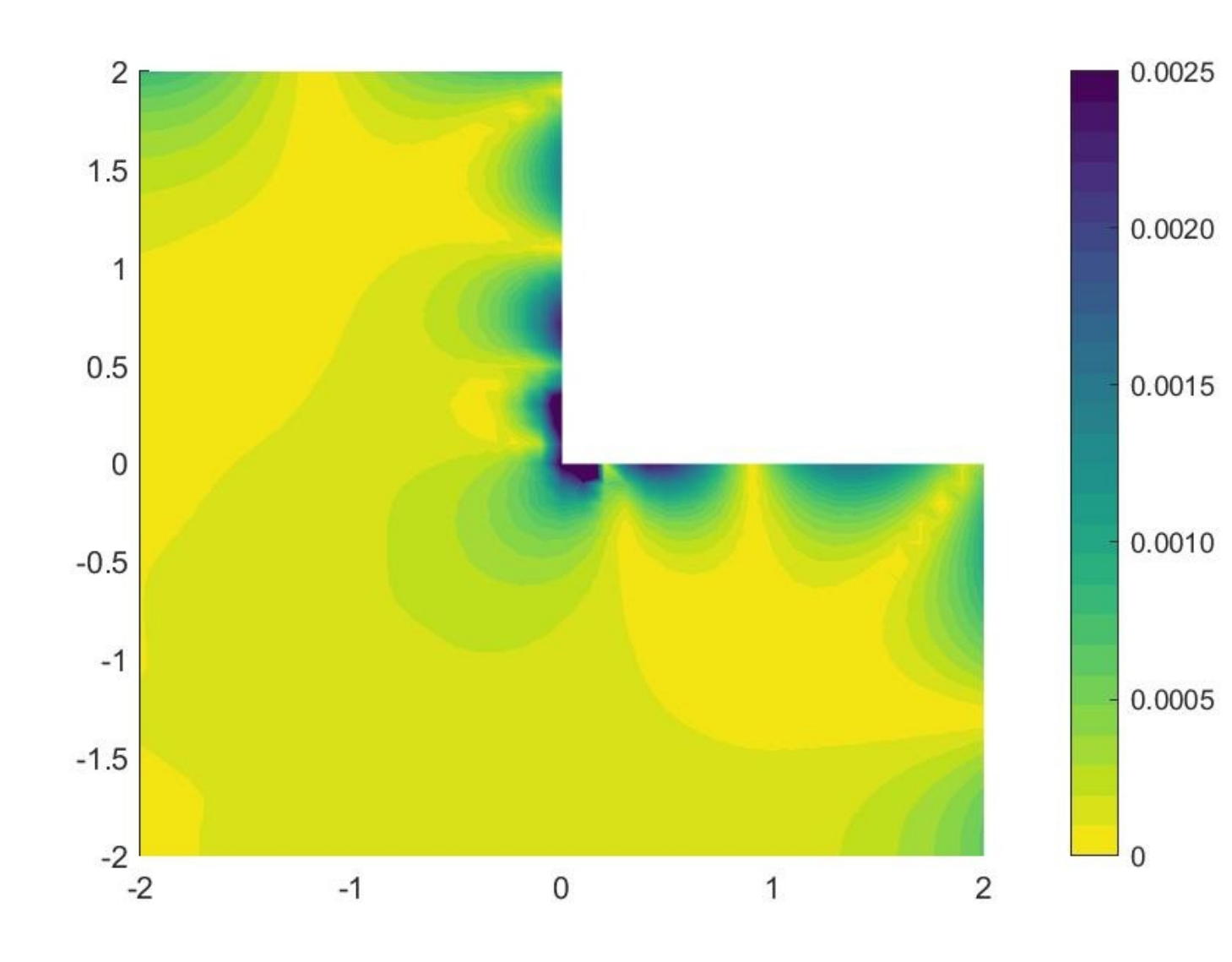}}
  \subfigure[PFNN-2]
  {\includegraphics[width=0.24\textwidth]
  {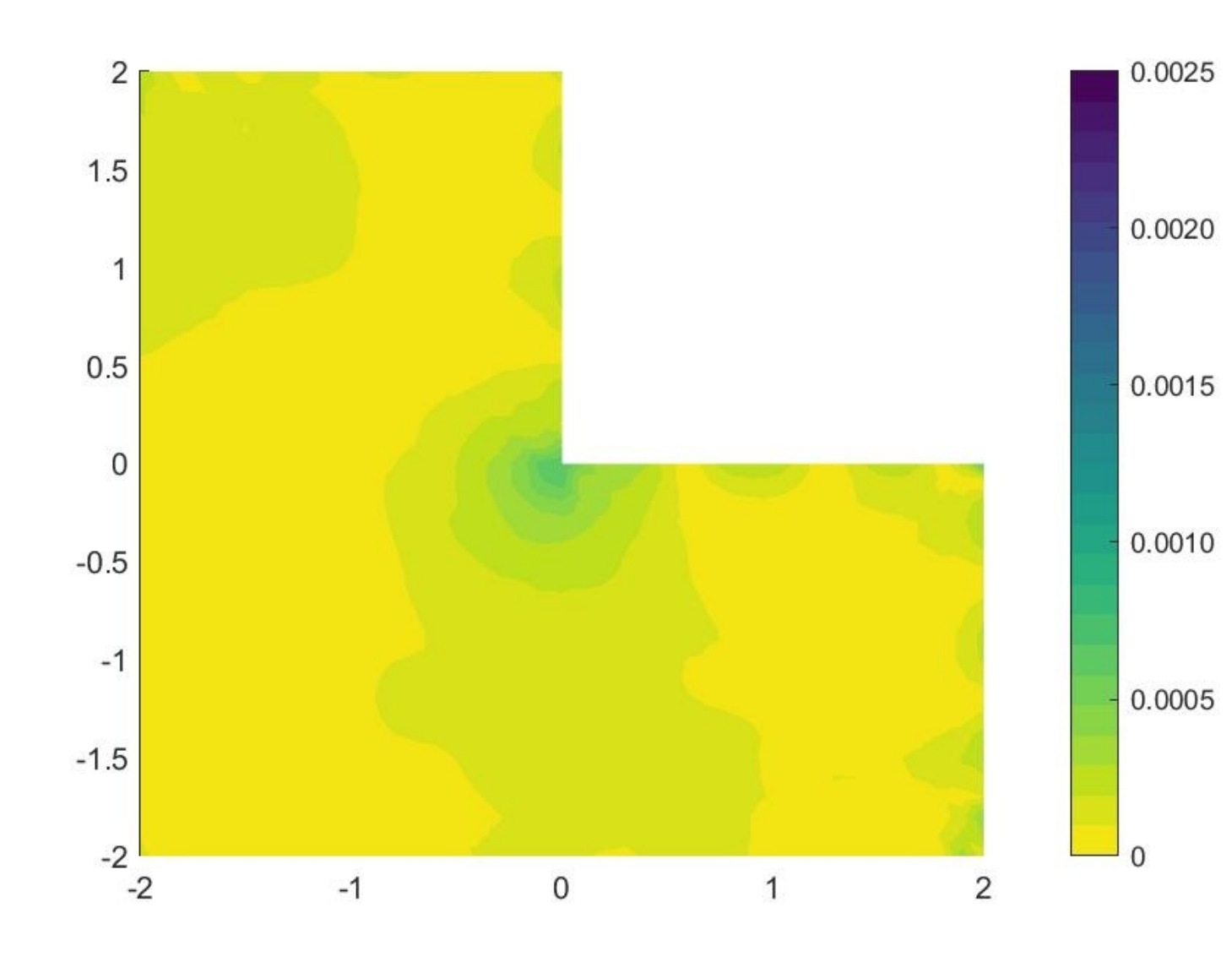}}
  \caption{Error distribution of the approximate solutions obtained with various methods for solving the anisotropic convection-diffusion equation on an L-shaped domain.}
  \label{f_linear_acd_error_distribution}
\end{figure}


The relative errors of the approximate solutions obtained by the tested algorithms are reported in Table \ref{t_linear_acd_results}, in which the influence of penalty factor for the least-squares neural network method and VarNet is also investigated.
It can be seen from the table that the two penalty-based methods is inferior to the classical finite element method in terms of accuracy for this test case, and
their performance is seriously affected by the change of penalty factors.
In contrast, PFNN-2 is free of any penalty term and is not subject to such influence, thus can outperform all the methods in terms of sustained accuracy.
Further, we draw the distribution of the solution errors of the tested algorithms in Figure \ref{f_linear_acd_error_distribution}, from which we observe that
there are significant errors occurring near the corner of the L-shaped domain when the finite element method and the penalty-based neural network approaches are used, indicating the inefficient treatment of corner singularity.
In addition, the two penalty-based neural network methods also have serious errors near the essential boundary.
PFNN-2, on the other hand, is free of such singularity and boundary effects, and can sustain the most accurate result among the tested methods, with the error uniformly distributed both within the domain and on the boundary.

\subsection{Allen-Cahn equation on a periodic domain}
The next test case is the well-known Allen-Cahn equation \cite{church2019high} for phase field simulation:
\begin{equation}
  \label{pde_ac}
  \dfrac{\partial u}{\partial t}
  -\epsilon^2 \Delta u +
  u^3 - u = 0,
\end{equation}
defined on a periodic domain $[0,2\pi]^2$.
We set the width of the transition layer to be $\epsilon=0.2$ and set the initial condition to be
\begin{equation*}
  u(\bm{x},0) = \tanh
  \dfrac{\sqrt{(x_1-\pi)^2 + (x_2-\pi)^2} - 2}
  {\epsilon\sqrt{2}}.
\end{equation*}
For this time-dependent problem, we conduct the simulation from $t=0$ to $t=10$.

The true solution of the above problem has some sharp gradients on the transition layer.
For dealing with the problem, domain decomposition based methods tend to be effective, since each network only need to learn one part of the solution.
To examine the effectiveness of the domain decomposition
strategies of various neural network methods including cPINNs, XPINNs, DeepDDM, D3M and PFNN-2,
four groups of experiments are conducted, in which the computational domain is decomposed into 1$\times$1, 3$\times$3, 5$\times$5 and 7$\times$7 sub-domains, respectively.
Configurations of the network structures for the various tested methods are listed in Table \ref{t_allen_cahn_networks}, where we try to maintain a similar amount of learnable parameters for each group of the test.
We sample 16000 points on the whole computational domain to form the test function sets.
The numbers of outer and inner iterations are set to be 150 and 40, respectively, resulting a total of 6000 iterations.
A numerical solution obtained by the Fourier spectral method with 1024 modes and second-order explicit Runge-Kutta temporal integrator with time step $\Delta t=0.001$ is taken as the ground truth for evaluating the solution error.
Different penalty coefficients including 1, 10, 100 and 1000 are examined for penalty-based methods and the value with the optimal performance is picked out for further analysis.

\begin{table}[!htb]
  \caption{Network structures of various methods for solving the Allen-Cahn equation.}
  \label{t_allen_cahn_networks}
  \centering
  \scalebox{0.78}{
  \renewcommand{\arraystretch}{1.2}
  \begin{tabu}{c|c c c}
    \tabucline[1pt]{-}
    $m$ & cPINNs/XPINNs/DeepDDM & D3M & PFNN-2 \\
    \hline
    \multirow{2}*{1}  & 1$\times$4 blocks (width:32) & 1$\times$4 blocks (width:32) & 1$\times$2$\times$4 blocks (width:22) \\
       & 7553 {DOFs} & 7619 DOFs & 7306 DOFs \\
    \hline
    \multirow{2}*{9}  & 9$\times$2 blocks (width:15) & 9$\times$2 blocks (width:15) & 9$\times$2$\times$2 blocks (width:10) \\
       & 7164 DOFs & 7452 DOFs & 6858 DOFs \\
    \hline
    \multirow{2}*{25} & 25$\times$1 blocks (width:14) & 25$\times$1 blocks (width:14) & 25$\times$2$\times$1 blocks (width:9) \\
       & 7025 DOFs & 7775 DOFs & 6800 DOFs \\
    \hline
    \multirow{2}*{49} & 49$\times$1 blocks (width:10) & 49$\times$1 blocks (width:10) & 49$\times$2$\times$1 blocks (width:6) \\
       & 7889 DOFs & 8967 DOFs & 7154 DOFs \\
    \tabucline[1pt]{-}
  \end{tabu}
  }
\end{table}

\begin{figure}[!htb]
  \centering
  \includegraphics[width=0.8\textwidth]
  {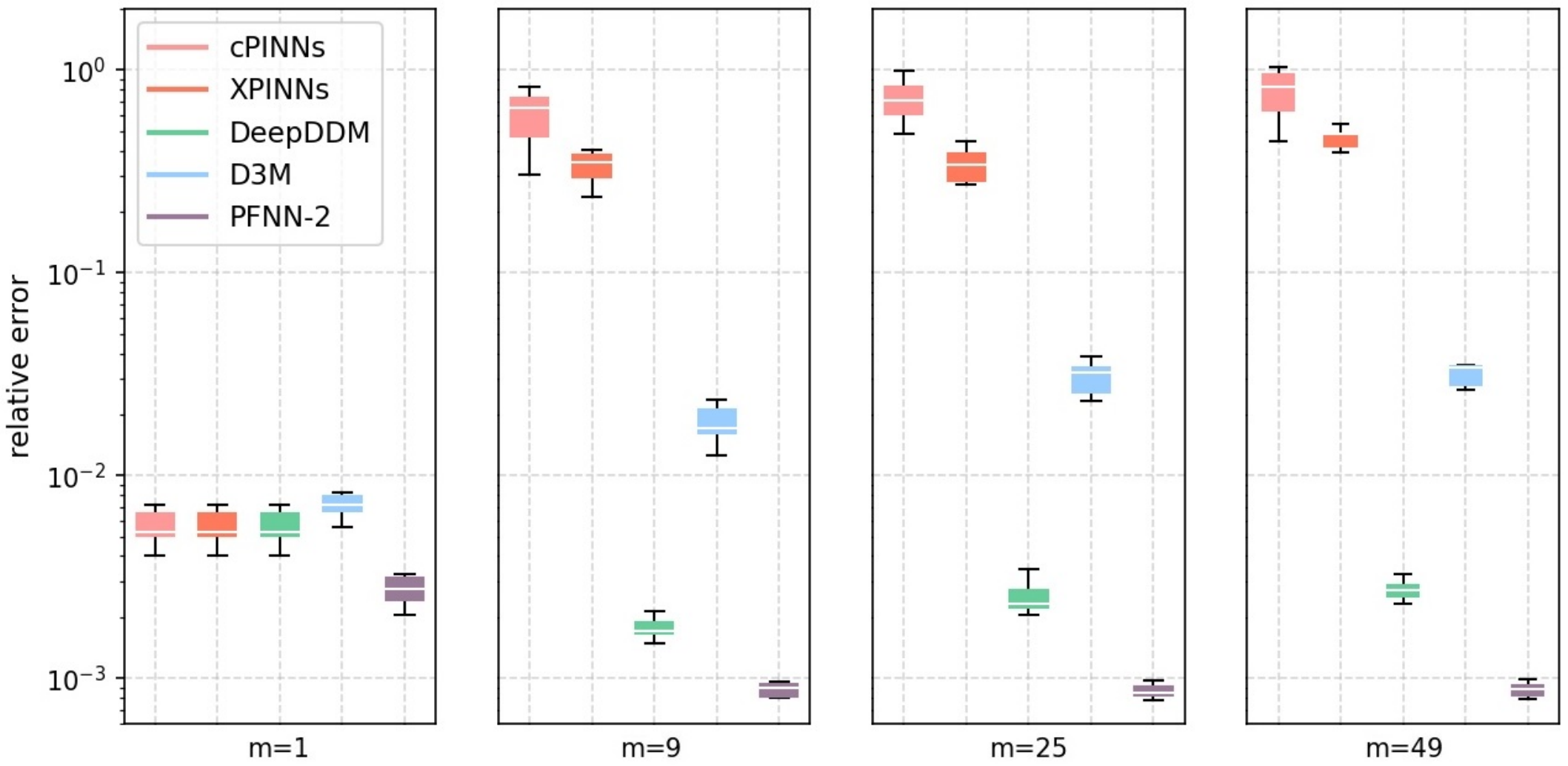}
  \caption{Relative errors of various domain decomposition based neural network methods for solving the Allen-Cahn equation.}
  \label{f_allen_cahn_accuracy}
\end{figure}

Figure \ref{f_allen_cahn_accuracy} provides an intuitive display of the numerical accuracy of different tested methods, in which box-plots are drawn to describe the five-number summary including the minimum observation, first quartile, middle value, third quartile, and maximum observation for each set of tests.
It can be seen that not all domain decomposition strategies can improve the accuracy. As the number of sub-domains increases, the achievable accuracy of most of the existing domain decomposition based neural network methods except DeepDDM gets worse, which reveals their drawbacks in handling multiple sub-problems with constraints.
Although DeepDDM can maintain a relatively satisfying accuracy, PFNN-2 has a more preferable performance, primarily due to its capability of dealing with smoothness constraints and boundary conditions. 

\subsection{Nonlinear anisotropic convection-diffusion equation on an irregular domain}
In this test case, we consider a nonlinear anisotropic convection-diffusion equation
\begin{equation}
  \label{pde_nacd}
  \dfrac{\partial u}{\partial t} -
  \nabla \cdot (\mathcal{A}\nabla u - \mathcal{B}) +
  \mathcal{C} = 0,\quad t\in[0,1]
\end{equation}
defined on an irregular domain with the shape of the Antarctica, as demonstrated in Table \ref{t_nonlinear_acd_network_and_domain}.
We set the exact solution to be {$u = \cos(\pi x_1) \cos(\pi x_2) \exp(-\pi t(x_1+x_2))$} that satisfies the Dirichlet boundary condition, and set the coefficients to be
\begin{equation*}
  \mathcal{A} = \left[
  \begin{array}{c c}
    u^2+1 & 2u^2 \\
    -2u^2 & u^2+1
  \end{array}
  \right],
  \quad
  \mathcal{B} = \left[
  \begin{array}{c}
    u^2 + 2u \\[0.25mm]
    -2u^2 + u
  \end{array}
  \right],
\end{equation*}
and $\mathcal{C}$ is derived according to (\ref{pde_nacd}).

\begin{table}[!htb]
  \caption{Network structures of various neural network methods 
  for solving the nonlinear anisotropic convection-diffusion equation on an irregular domain  with domain decomposition.}
  \label{t_nonlinear_acd_network_and_domain}
  \centering
  \scalebox{0.74}{
  \renewcommand{\arraystretch}{1.3}
  \begin{tabu}{c|c c c c}
    \tabucline[1pt]{-}
    $m$ & 1 & 14 & 26 & 42 \\
    \hline
    Network (cPINNs/ & 1$\times$4 blocks (width:30) & 14$\times$2 blocks (width:11) & 26$\times$1 blocks (width:13) & 42$\times$1 blocks (width:10) \\
    DeepDDM) & 6661 DOFs & 6328 DOFs & 6448 DOFs & 6762 DOFs \\
    \hline
    \multirow{2}*{Network (PFNN-2)} & 1$\times$2$\times$4 blocks (width: & 14$\times$2$\times$2 blocks (width: & 26$\times$2$\times$1 blocks (width: & 42$\times$2$\times$1 blocks (width: \\
    & 20) 6082 DOFs & 7) 5712 DOFs & 8) 5876 DOFs & 6) 6132 DOFs \\
    \hline
    \multirow{6}*{\shortstack{Decomposition \\ of the domain}}
    & \multirow{6}*{\includegraphics[width=0.25\textwidth]{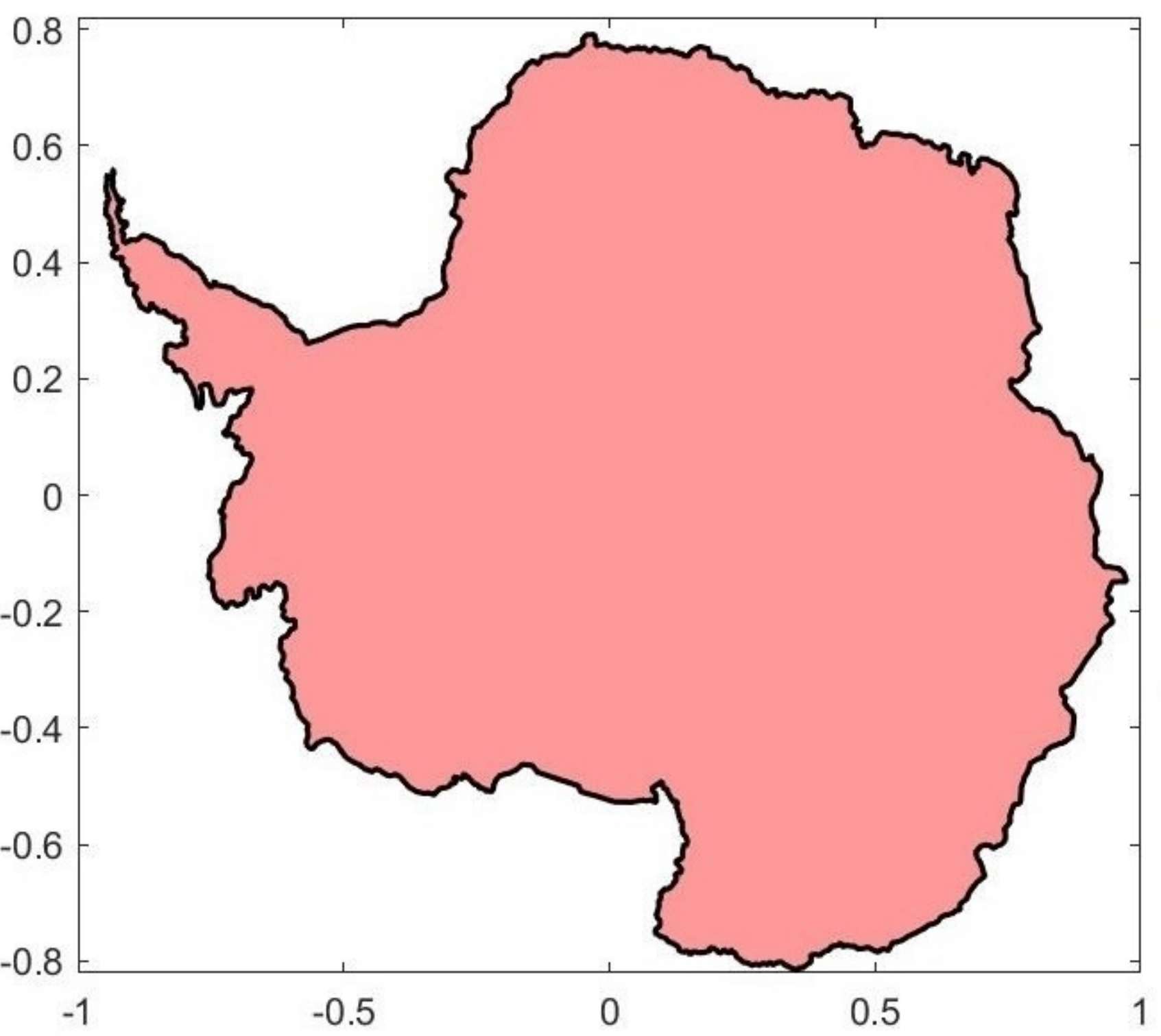}}
    & \multirow{6}*{\includegraphics[width=0.25\textwidth]{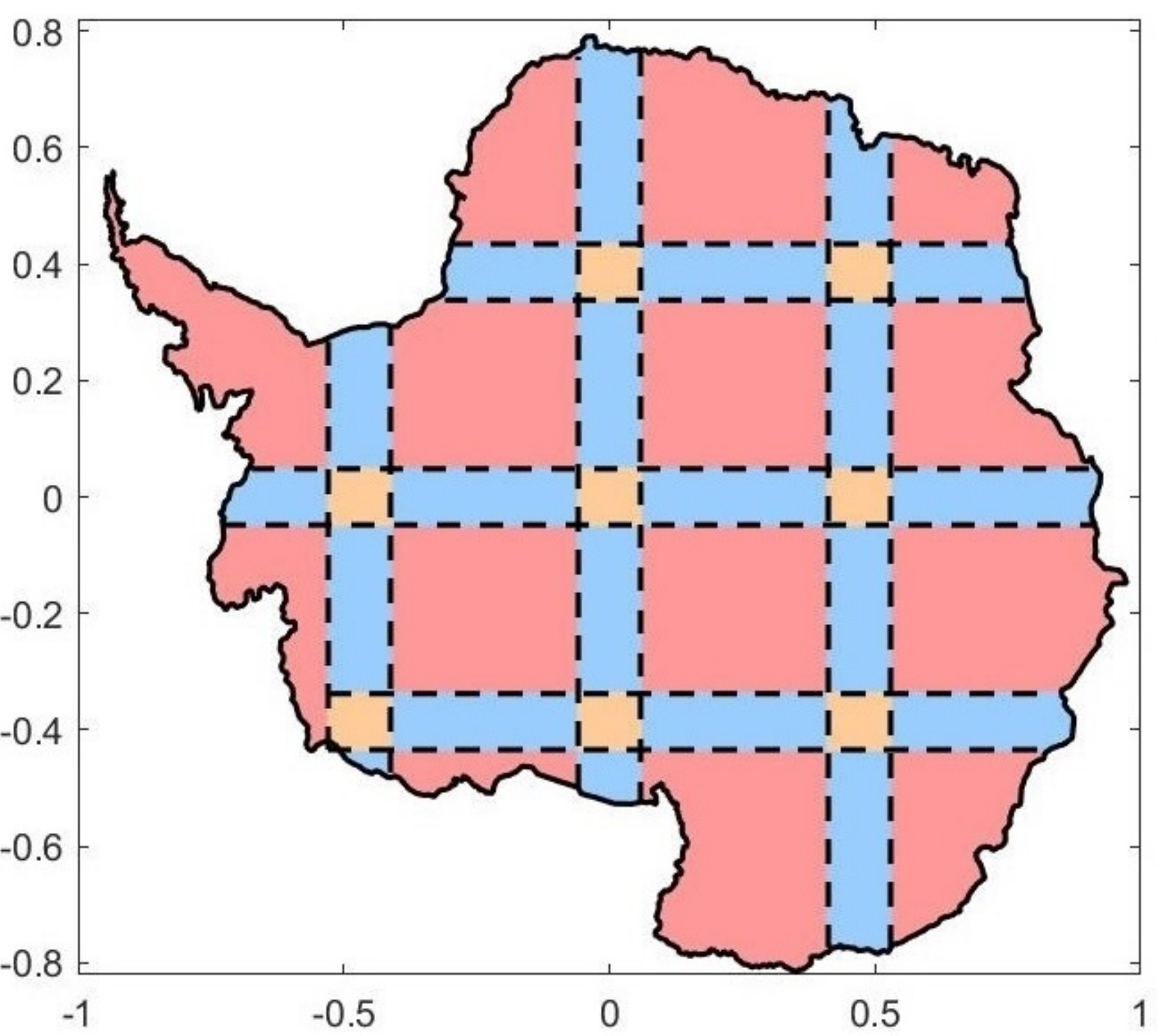}}
    & \multirow{6}*{\includegraphics[width=0.25\textwidth]{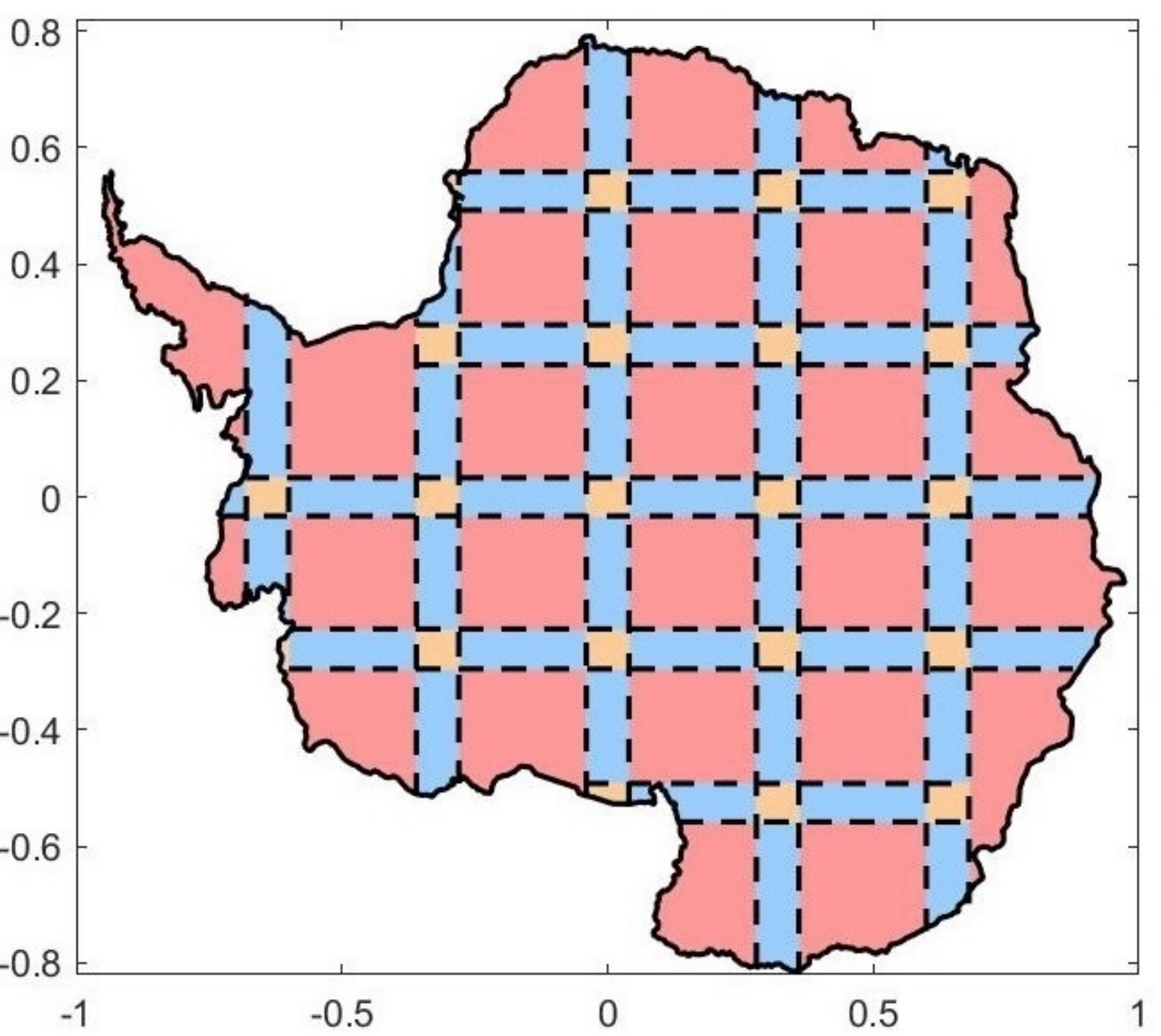}}
    & \multirow{6}*{\includegraphics[width=0.25\textwidth]{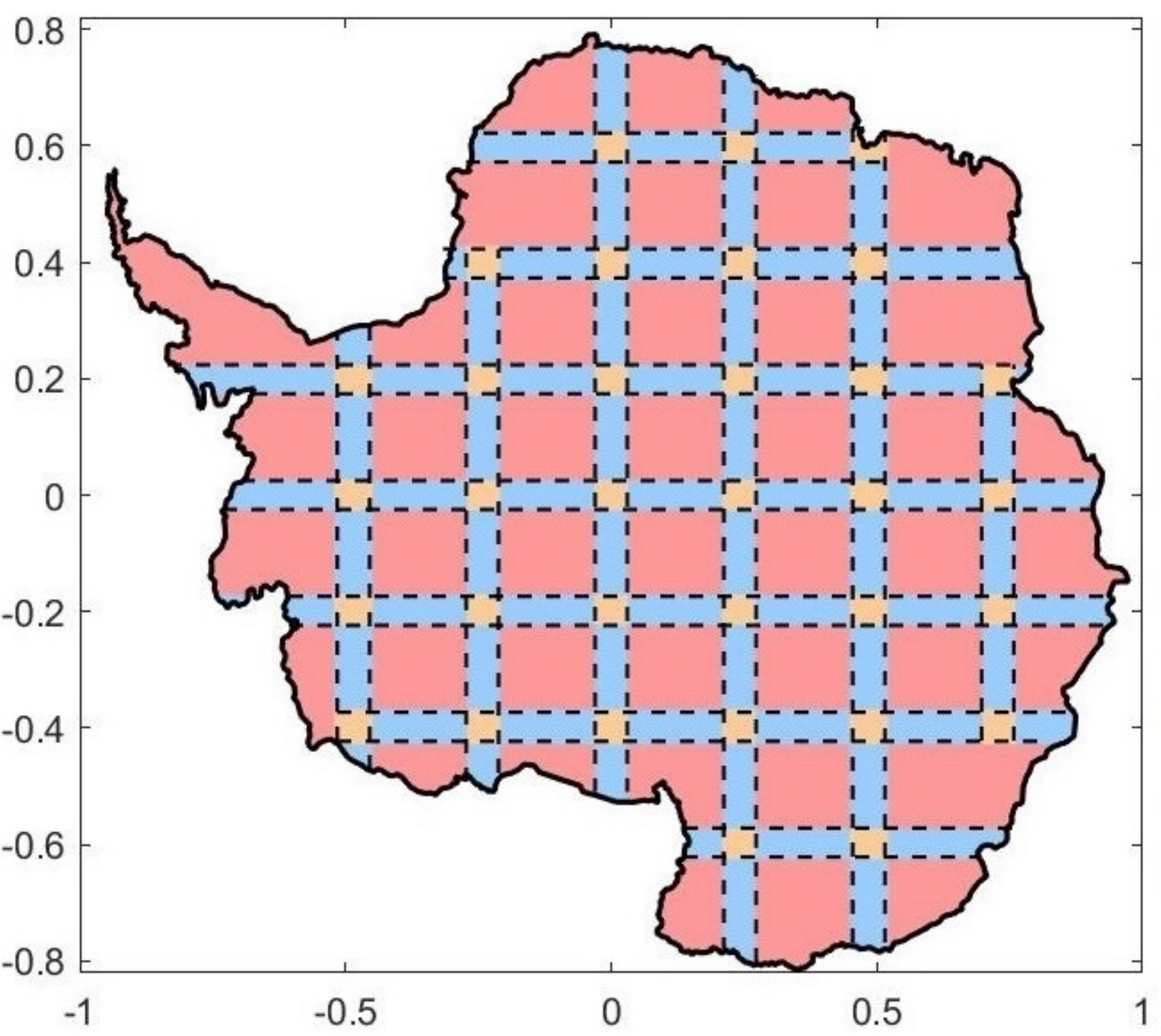}}
    \\ & \\ & \\ & \\ & \\ & \\
    \tabucline[1pt]{-}
  \end{tabu}
  }
\end{table}


To investigate the effectiveness of domain decomposition strategies of different neural network methods, we examine the performance of PFNN-2 and compare it with cPINNs and DeepDDM, which are selected as representative non-overlapping and overlapping domain decomposition neural network methods, respectively.
We conduct four groups of experiments, in which the computational domain is decomposed into 1, 14, 26 and 42 sub-domains, respectively.
In each group of experiments, we set network structure of each method with similar amount of undecided parameters, as illustrated in Table \ref{t_nonlinear_acd_network_and_domain}.
A total of 6000 points are sampled on the whole computational domain to form the test function set.
Other experiment configurations are the same to  the previous experiment.

\begin{figure}[!htb]
  \centering
  \subfigure[cPINNs and DeepDDM]
  {\includegraphics[width=0.4\textwidth]
  {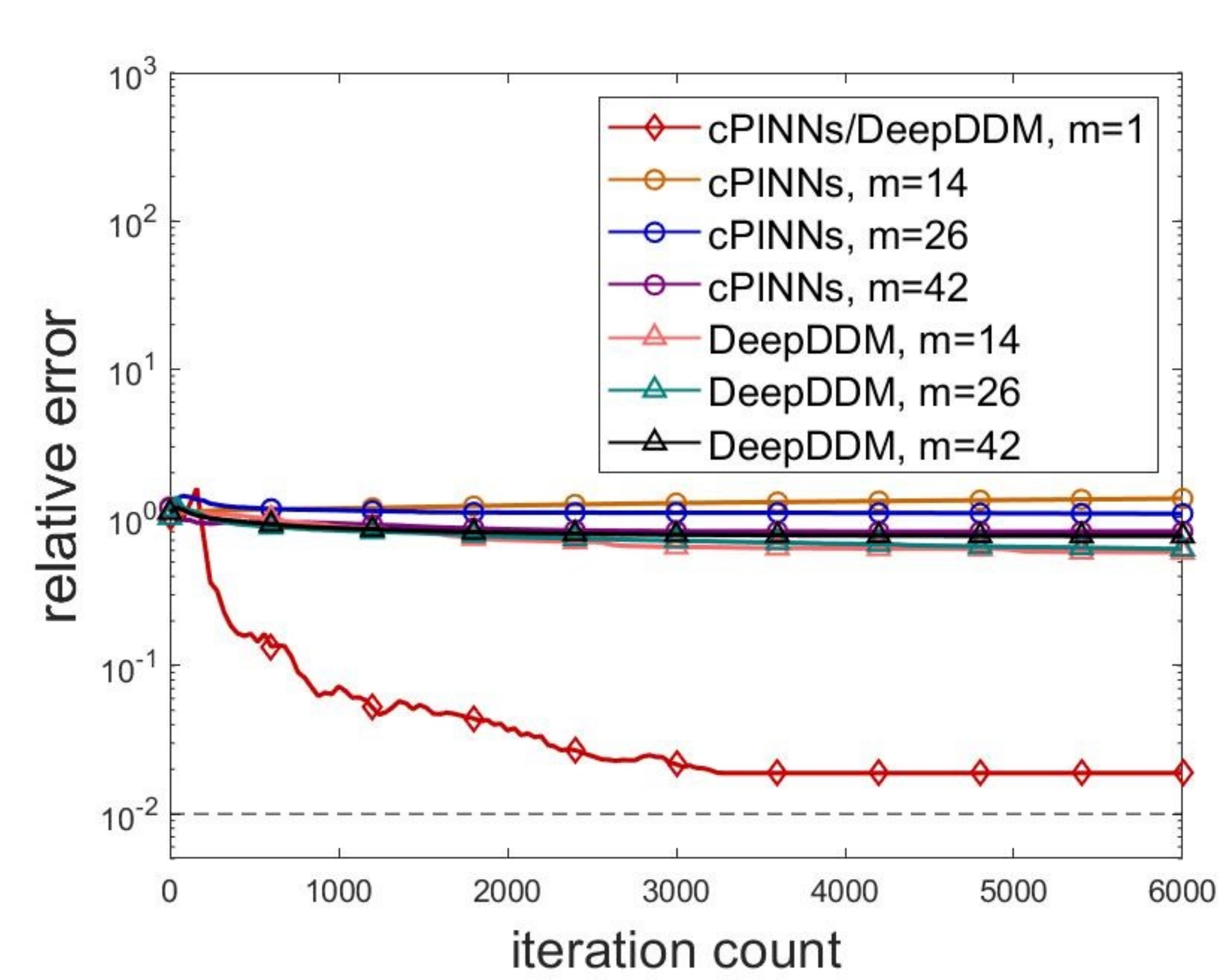}}
  \quad
  \subfigure[PFNN-2]
  {\includegraphics[width=0.4\textwidth]
  {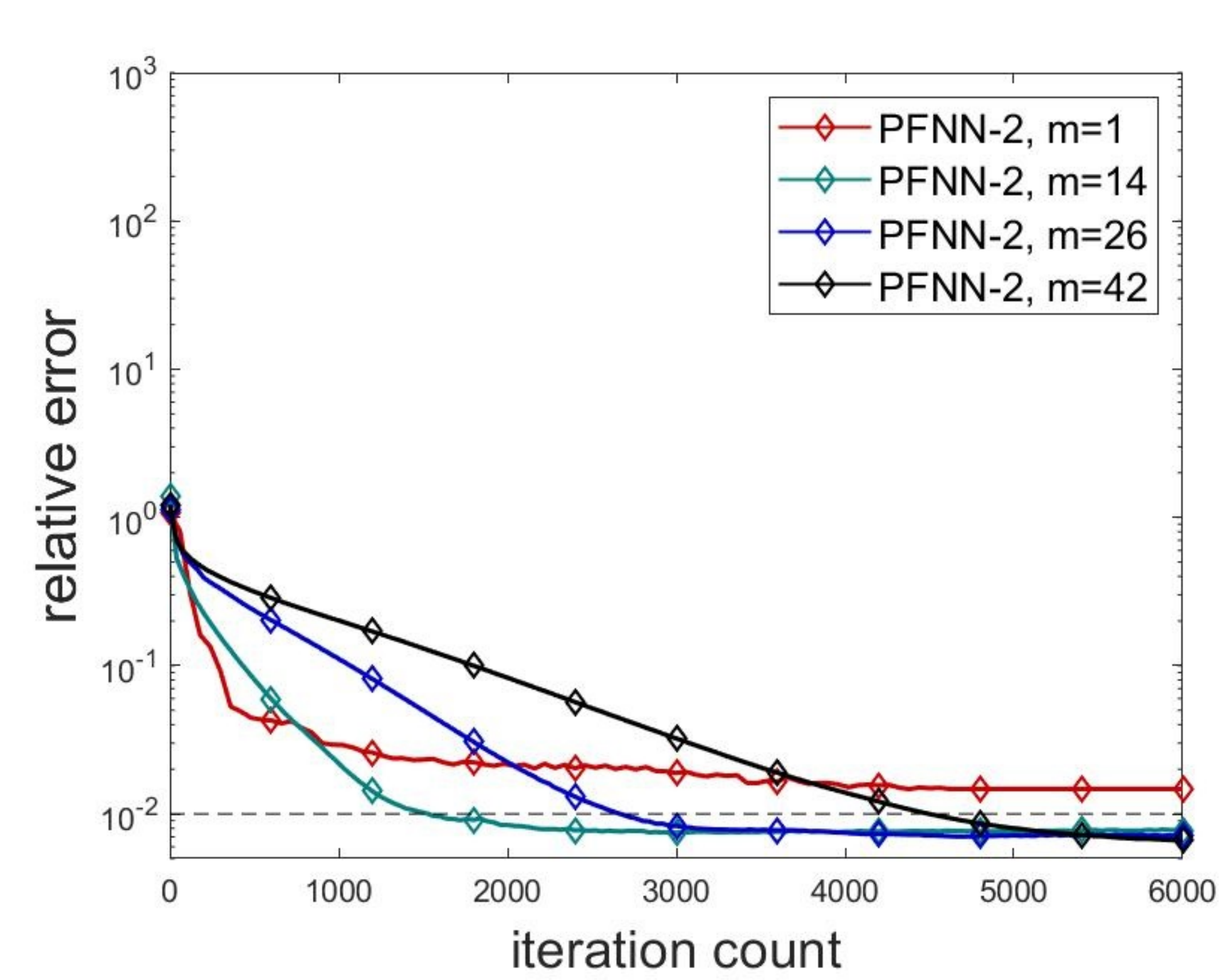}}
  \caption{Error evolution histories of various methods during the training process in solving the nonlinear anisotropic convection-diffusion equation on an Antarctica-shaped domain.}
  \label{f_nonlinear_acd_error_history}
\end{figure}

Figure \ref{f_nonlinear_acd_error_history} demonstrates the error evolution histories of various methods.
It can be seen from the figure that it is hard for cPINNs and DeepDDM to converge to the true solution in the cases that multiple sub-domains are employed, which exposes their defects in tackling problem with strong nonlinearity.
By comparison, thanks to the avoidance of penalty terms, as well as the reduction of the smoothness constraints and the flexble treatment of the complex geometry, PFNN-2 is more suitable for solving this problem and can achieve more accurate and stable results.

\subsection{Viscous Burgers equation on a 3D cube}
The last test case is a nonlinear equation system -- the viscous Burgers equation
\begin{equation}
  \label{PDE_Burgers}
  \dfrac{\partial \bm{u}}{\partial t}
  - \dfrac{1}{\mbox{Re}} \Delta \bm{u}
  + \bm{u}\cdot\nabla\bm{u}
  = \bm{0},
\end{equation}
defined on domain $[0,1]^3$ with Dirichlet boundary condition defined on $\Gamma_D = [0,1]^2\times\{0,1\}$ and Neumann boundary condition defined on $\Gamma_N = \Gamma\backslash\Gamma_D$, where the Reynolds number is $\mbox{Re}=100$.
We set the true solution to be
\begin{equation*}
  \bm{u} =
  \left[
  \begin{array}{l}
    u_1(\bm{x},t) \\[3.5mm]
    u_2(\bm{x},t) \\[3.5mm]
    u_3(\bm{x},t)
  \end{array}
  \right]
  =
  \left[
  \scalebox{0.78}{$
  \begin{array}{l}
    -\dfrac{2}{\mbox{Re}}
    \left( \dfrac{1 + \pi\cos(\pi x_1) \sin(\pi x_2) \sin(\pi x_3) \exp(-3\pi^2 t/\mbox{Re})}
           {1 + x_1 + \sin(\pi x_1) \sin(\pi x_2) \sin(\pi x_3) \exp(-3\pi^2 t/\mbox{Re})} \right) \\[5.5mm]
    -\dfrac{2}{\mbox{Re}}
    \left( \dfrac{\pi\sin(\pi x_1) \cos(\pi x_2) \sin(\pi x_3) \exp(-3\pi^2 t/\mbox{Re})}
           {1 + x_1 + \sin(\pi x_1) \sin(\pi x_2) \sin(\pi x_3) \exp(-3\pi^2 t/\mbox{Re})} \right) \\[5.5mm]
    -\dfrac{2}{\mbox{Re}}
    \left( \dfrac{\pi\sin(\pi x_1) \sin(\pi x_2) \cos(\pi x_3) \exp(-3\pi^2 t/\mbox{Re})}
           {1 + x_1 + \sin(\pi x_1) \sin(\pi x_2) \sin(\pi x_3) \exp(-3\pi^2 t/\mbox{Re})} \right)
  \end{array}
  $}
  \right],
\end{equation*}
and conduct the simulation from $t=0$ to $t=1.0$.

\begin{figure}[!htb]
  \centering
  \includegraphics[width=0.75\textwidth]
  {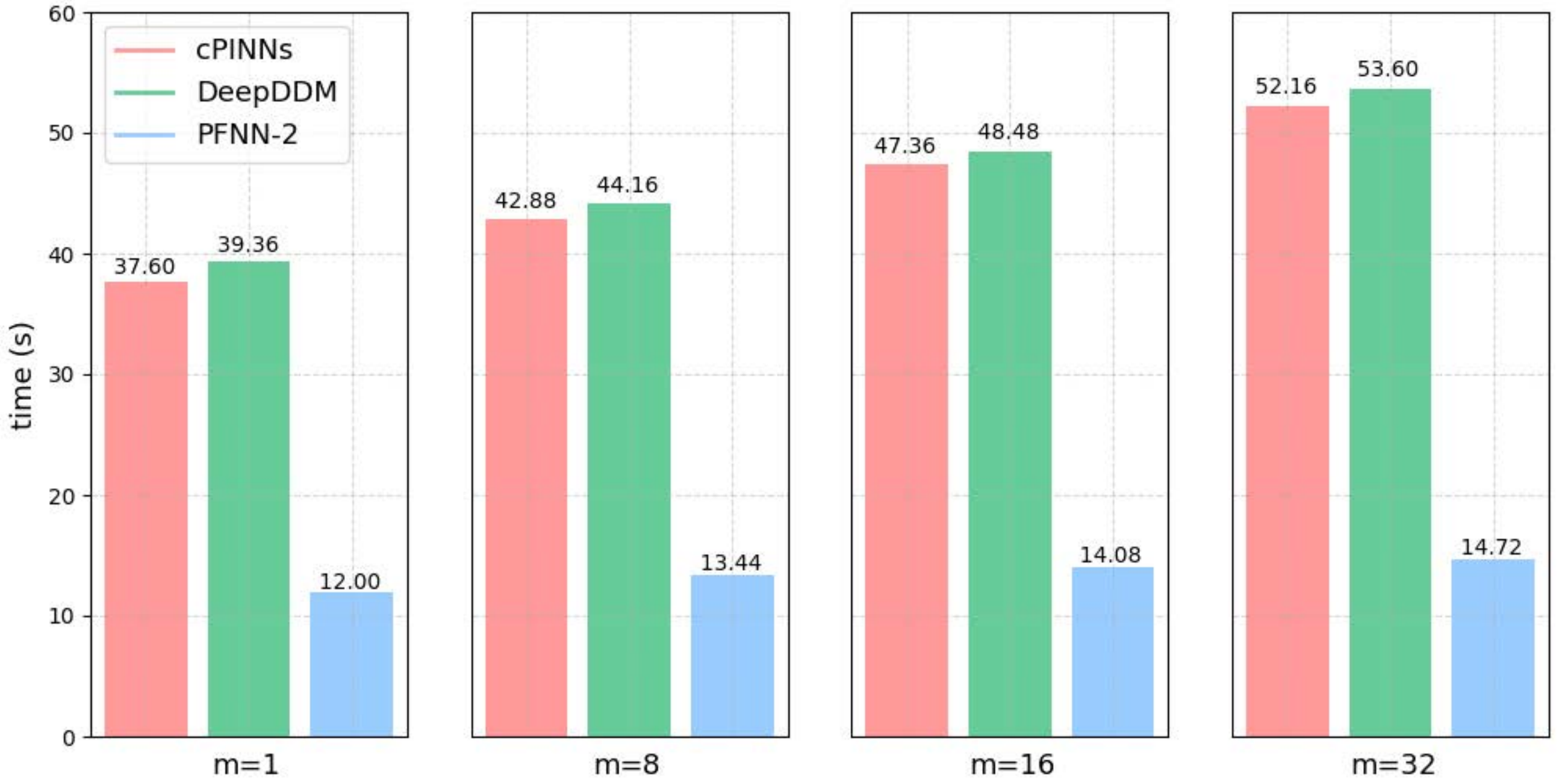}
  \caption{Efficiency of various domain decomposition based neural network methods in solving the viscous Burgers equations in weak scaling efficiency test.}
  \label{f_burgers_efficiency_weak}
\end{figure}

In this test case, we examine both accuracy and efficiency of various methods including cPINNs, DeepDDM and PFNN-2.
Firstly, the weak scaling efficiency of the tested methods is investigated.
We conduct four sets of experiments, in which the computational domain is decomposed into 1, $2\times2\times2$, $4\times2\times2$ and $4\times4\times2$ sub-domains, respectively.
We keep the scale of task assigned to each processor similar.
Specifically, we sample 5120 points on each sub-domain respectively to form the test function set, and adopt networks with similar scale to approximate local solution on each sub-domain, as illustrated in Table \ref{t_burgers_accuracy_weak}.
Other experiment configurations are the same to the previous experiment.
The accuracy of various method is reported in Table \ref{t_burgers_accuracy_weak}.
Not surprisingly, the accuracy of all the methods is improved with the increase of the number of sub-domains, since the representative capability of approximate solution also grows stronger.
In particular, PFNN-2 acquires the most accurate approximate solutions in all the test sets.
Moveover, Figure \ref{f_burgers_efficiency_weak} demonstrates the runtime for 400 iterations, which indicates that all the tested algorithms can maintain a relatively stable efficiency. In all the experiment groups, PFNN-2 always achieves the highest efficiency.

\begin{table}[!htb]
  \caption{Accuracy of various domain decomposition based neural network methods for solving the viscous Burgers equations in weak scaling efficiency test.}
  \label{t_burgers_accuracy_weak}
  \centering
  \scriptsize
  \renewcommand{\arraystretch}{1.2}
  \begin{tabu}{c|c c c c c}
    \tabucline[1pt]{-}
    Method & Network structure & DOFs & Errors ($u_1$) & Errors ($u_2$) & Errors ($u_3$) \\
    \hline
    \multirow{4}*{cPINNs}  & 1$\times$2 blocks (width:15)  & 1$\times$828  & 7.31e-02$\pm$2.57e-02 & 1.31e-01$\pm$3.64e-02 & 1.34e-01$\pm$4.08e-02 \\
                           & 8$\times$2 blocks (width:15)  & 8$\times$828  & 3.09e-02$\pm$8.38e-03 & 4.54e-02$\pm$1.33e-02 & 5.39e-02$\pm$1.63e-02 \\
                           & 16$\times$2 blocks (width:15) & 16$\times$828 & 2.64e-02$\pm$7.40e-03 & 3.92e-02$\pm$1.02e-02 & 4.58e-02$\pm$1.36e-02 \\
                           & 32$\times$2 blocks (width:15) & 32$\times$828 & 2.33e-02$\pm$6.92e-03 & 3.31e-02$\pm$8.19e-03 & 3.72e-02$\pm$8.76e-03 \\
    \hline
    \multirow{4}*{DeepDDM} & 1$\times$2 blocks (width:15)  & 1$\times$828  & 7.31e-02$\pm$2.57e-02 & 1.31e-01$\pm$3.64e-02 & 1.34e-01$\pm$4.08e-02 \\
                           & 8$\times$2 blocks (width:15)  & 8$\times$828  & 3.38e-02$\pm$1.15e-02 & 5.55e-02$\pm$1.58e-02 & 6.45e-02$\pm$2.12e-02 \\
                           & 16$\times$2 blocks (width:15) & 16$\times$828 & 2.46e-02$\pm$6.74e-03 & 3.68e-02$\pm$1.13e-02 & 4.74e-02$\pm$1.48e-02 \\
                           & 32$\times$2 blocks (width:15) & 32$\times$828 & 1.64e-02$\pm$5.41e-03 & 2.38e-02$\pm$7.20e-03 & 2.95e-02$\pm$7.87e-03 \\
    \hline
    \multirow{4}*{PFNN-2}  & 1$\times$2$\times$2 blocks (width:10)  & 1$\times$806  & 5.06e-02$\pm$1.36e-02 & 8.02e-02$\pm$2.14e-02 & 7.67e-02$\pm$2.05e-02 \\
                           & 8$\times$2$\times$2 blocks (width:10)  & 8$\times$806  & 2.23e-02$\pm$4.94e-03 & 2.76e-02$\pm$5.46e-03 & 3.22e-02$\pm$6.44e-03 \\
                           & 16$\times$2$\times$2 blocks (width:10) & 16$\times$806 & 9.54e-03$\pm$2.27e-03 & 1.43e-02$\pm$3.69e-03 & 1.51e-02$\pm$3.83e-03 \\
                           & 32$\times$2$\times$2 blocks (width:10) & 32$\times$806 & \textbf{6.72e-03$\pm$1.80e-03} & \textbf{9.72e-03$\pm$2.58e-03} & \textbf{9.81e-03$\pm$2.28e-03} \\
    \tabucline[1pt]{-}
  \end{tabu}
\end{table}

\begin{figure}[!htb]
  \centering
  \includegraphics[width=0.75\textwidth]
  {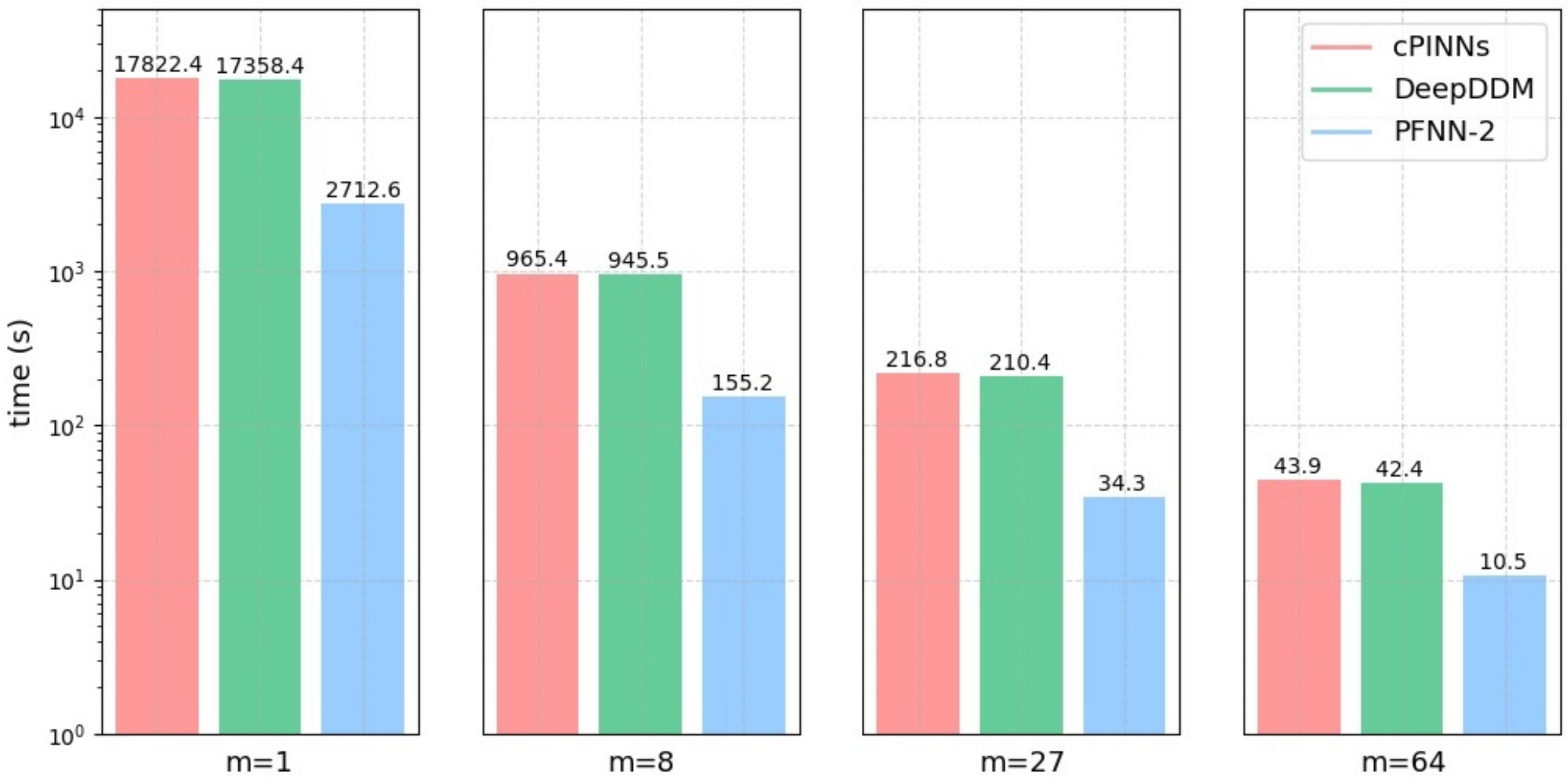}
  \caption{Efficiency of various domain decomposition based neural network methods in solving the viscous Burgers equations on a 3D cube.}
  \label{f_burgers_efficiency_strong}
\end{figure}

Next, the strong scaling efficiency of the tested approaches is examined.
Four sets of experiments are carried out, in which the computational domain is split into $1^3$, $2^3$, $3^3$ and $4^3$ sub-domains, respectively.
In this test, the overall scale of problem is kept similar in different test groups.
To be specific, we sample 160,000 points on the whole computational domain to form the test function set,
and keep the amounts of learnable parameters of different approaches in each group in a similar range, as shown in Table \ref{t_burgers_accuracy_strong}.
Other experiment configurations are the same to the previous experiment.
Table \ref{t_burgers_accuracy_strong} also reports the solution errors of the various tested algorithms, which once more demonstrates the advantage of PFNN-2 with respect to the numerical accuracy.
To further compare the performance in terms of computing time, we draw the runtime cost for 400 iterations in Figure \ref{f_burgers_efficiency_strong}, from which it can be seen that with the increase of the number of processors, the time-to-solution of all the methods is enhanced significantly.
Among all the tested algorithms, PFNN-2 is always the most efficient.

\begin{table}[!htb]
  \caption{Accuracy of various domain decomposition based neural network methods for solving the viscous Burgers equations in strong scaling efficiency test.}
  \label{t_burgers_accuracy_strong}
  \centering
  \scriptsize
  \renewcommand{\arraystretch}{1.2}
  \begin{tabu}{c|c c c c c}
    \tabucline[1pt]{-}
    Method & Network structure & DOFs & Errors ($u_1$) & Errors ($u_2$) & Errors ($u_3$) \\
    \hline
    \multirow{4}*{cPINNs}  & 1$\times$5 blocks (width:72)  & 47883 & 6.51e-02$\pm$2.63e-02 & 1.21e-01$\pm$3.77e-02 & 1.25e-01$\pm$3.85e-02 \\
                           & 8$\times$3 blocks (width:33)  & 47016 & 3.36e-02$\pm$8.34e-03 & 4.78e-02$\pm$1.19e-02 & 5.39e-02$\pm$1.37e-02 \\
                           & 27$\times$2 blocks (width:22) & 45819 & 2.89e-02$\pm$6.73e-03 & 3.77e-02$\pm$8.84e-03 & 4.05e-02$\pm$1.05e-02 \\
                           & 64$\times$2 blocks (width:14) & 47680 & 2.70e-02$\pm$6.48e-03 & 4.03e-02$\pm$9.32e-03 & 3.92e-02$\pm$8.39e-03 \\
    \hline
    \multirow{4}*{DeepDDM} & 1$\times$5 blocks (width:72)  & 47883 & 6.51e-02$\pm$2.63e-02 & 1.21e-01$\pm$3.77e-02 & 1.25e-01$\pm$3.85e-02 \\
                           & 8$\times$3 blocks (width:33)  & 47016 & 3.81e-02$\pm$7.81e-03 & 4.76e-02$\pm$1.25e-02 & 6.13e-02$\pm$1.42e-02 \\
                           & 27$\times$2 blocks (width:22) & 45819 & 1.60e-02$\pm$4.25e-03 & 2.31e-02$\pm$5.30e-03 & 2.58e-02$\pm$5.22e-03 \\
                           & 64$\times$2 blocks (width:14) & 47680 & 1.21e-02$\pm$3.26e-03 & 1.75e-02$\pm$4.82e-03 & 1.72e-02$\pm$3.87e-03 \\
    \hline
    \multirow{4}*{PFNN-2}  & 1$\times$2$\times$5 blocks (width:50)  & 46706 & 4.59e-02$\pm$1.30e-02 & 6.73e-02$\pm$1.46e-02 & 5.92e-02$\pm$1.32e-02 \\
                           & 8$\times$2$\times$3 blocks (width:23)  & 47152 & 1.15e-02$\pm$3.26e-03 & 1.72e-02$\pm$4.56e-03 & 1.67e-02$\pm$4.30e-03 \\
                           & 27$\times$2$\times$2 blocks (width:15) & 45522 & 6.49e-03$\pm$1.94e-03 & 9.37e-03$\pm$2.86e-03 & 9.22e-03$\pm$1.77e-03 \\
                           & 64$\times$2$\times$1 blocks (width:15) & 46464 & \textbf{5.12e-03$\pm$1.63e-03} & \textbf{7.84e-03$\pm$2.03e-03} & \textbf{7.68e-03$\pm$1.96e-03} \\
    \tabucline[1pt]{-}
  \end{tabu}
\end{table}

\section{Conclusion}
In this paper, a new penalty-free neural network method --- PFNN-2 is presented for solving PDEs including non-self-adjoint time-dependent ones.
PFNN-2 employs a Galerkin variational principle to convert the original equation into a weak form with low smoothness requirement, splits the problem with complex initial-boundary constraints into two simpler unconstrained ones, and adopts two neural networks to approximate the corresponding local solutions, which are then combined to form the approximate solution on the whole domain.
A domain decomposition strategy is introduced in PFNN-2 to accelerate the training efficiency while maintaining the numerical accuracy.
Experiment results demonstrate that PFNN-2 can surpass existing state-of-the-art in terms of accuracy, efficiency and scalability.
Possible future explorations on PFNN-2 would involve the improvement on the network structure, the investigation on the error estimation, and the study on solving more challenging problems such as the ones with strong singularities.

\section*{Acknowledgments}

This work was supported in part by National Natural Science Foundation of China (grant\# 12131002).


\end{document}